\documentclass[english]{article}
\usepackage[T1]{fontenc}
\usepackage[latin9]{inputenc}
\usepackage{textcomp}
\usepackage{ifsym}
\usepackage{slashed}
\usepackage{amsmath}
\usepackage{amssymb}

\makeatletter

\newcommand{\lyxmathsym}[1]{\ifmmode\begingroup\def\b@ld{bold}
  \text{\ifx\math@version\b@ld\bfseries\fi#1}\endgroup\else#1\fi}

\newcommand{\lyxaddress}[1]{
\par {\raggedright #1
\vspace{1.4em}
\noindent\par}
}

\makeatother

\usepackage{babel}

\begin{document}

\title{UNIVERSAL PAD\'{E} APPROXIMATION}

\author{Nicholas J.  Daras$^{1}$ and Vassili Nestoridis$^{2}$}

\maketitle

\lyxaddress{$^{1}$Department of Mathematics, University of Military Education-
Hellenic Army Academy, Vari \textendash{} 16673, Greece, e-mail: darasn@sse.gr, }

\lyxaddress{$^{2}$Department of Mathematics, University of Athens, Panepistemiopolis,
15784, Athens, Greece, e-mail: vnestor@math.uoa.gr}
\begin{abstract}
In transferring some results from universal Taylor series to the case
of Pad\'{e} approximants we obtain stronger results, such as, universal
approximation on compact sets of arbitrary connectivity and generic
results on planar domains of any connectivity and not just on simply
connected domains.

\noindent \textbf{Subject Classification MSC2010}: primary 41A21,
30K05, secondary 30B10, 30E10, 30K99, 41A10, 41A20
\end{abstract}


\section{Introduction}

Professor C. Brezinski in a colloquium talk at the University of Athens
held on September 22, 2010, presented an overview on Pad\'{e} approximants.
In his thesis, under Charles Hermite, Henri Pad\'{e} arranged these approximants
in a double array now known as the Pad\'{e} table of the formal power
series $f={\displaystyle \sum_{v=0}^{\infty}a_{v}z^{n}=(a_{0},a_{1},...}):$\[
(\left[p/q\right]_{f})_{p,q=0,1,2,\ldots}.\]

In particular, for $q=0$ the Pad\'{e} approximant $\Bigl[p/0\Bigr]_{f}\;(p=0,1,2,...)$
is a polynomial and coincides with a partial sum of $f$ , while for
$q\geq1$ the Pad\'{e} approximant $\Bigl[p/q\Bigr]_{f}$ is a rational
function with some poles in general.

Pad\'{e} approximants have been applied in proofs of irrationality and
transcendence in number theory, in practical computation of special
functions, and in the analysis of different schemes for numerical
solution of ordinary or partial differential equations. For a short
history on Pad\'{e} approximation we refer to the review-research paper
\cite{r18} \cite{r18}.

In addition to their wide variety of applications, Pad\'{e} approximants
are also connected with continued fraction expansions (\cite{r29}, \cite{r30}, \cite{r46} and \cite{r51}),
orthogonal polynomials (\cite{r9}, \cite{r22}, \cite{r50} and \cite{r53}), moment
problems (\cite{r3} and cite{r44}), the theory of quadrature (\cite{r12} and \cite{r13})
and convergence acceleration methods (\cite{r10}, \cite{r11} and \cite{r53}). However,
the application which brought them to prominence in the 1960\textquoteright{}s
and 1970\textquoteright{}s, was localizing the singularities of functions:
in various problems, for example in inverse scattering theory, one
would have a means for computing the coefficients of a power series
$f$ . One could use these coefficients to compute a Pad\'{e} approximant
to $f$, and use the poles of the approximant as predictors of the
location of poles or other singularities of $f$ . Under certain conditions
on $f$, which were often satisfied in physical examples, this process
could be theoretically justified.

Furthermore, since 1965 a growing interest for Pad\'{e} approximants appeared
in theoretical physics, chemistry, electronics and numerical analysis;
see, for example, the books \cite{r2}, \cite{r4}, \cite{r13} and \cite{r26} and the international
conferences proceedings \cite{r20}, \cite{r48} and \cite{r54}.

One of the fascinating features of Pad\'{e} approximants is the complexity
of their asymptotic behaviour (\cite{r1}, \cite{r2}, \cite{r3}, \cite{r5}, \cite{r7}, \cite{r15}, \cite{r19}, \cite{r21}, \cite{r23}, \cite{r24}, \cite{r30}, \cite{r34}, \cite{r45}, \cite{r47}, \cite{r52} and \cite{r56}).
The convergence problem for Pad\'{e} approximants
can be stated as follows. Given a power series $f={\displaystyle \sum_{v=0}^{\infty}a_{v}z^{v}=\left(a_{0},a_{1},\cdots\right)}$
examine the convergence of subsequences $\Bigl(\Bigl[p_{n}/q_{n}\Bigr]_{f}\Bigr)_{n=0,1,2,...}$
extracted from the Pad\'{e} table as $n\rightarrow+\infty$. If $\mathcal{\mathcal{\mathcal{\mathtt{\mathrm{\,\mathfrak{f}}}}}}$
is the Taylor development of an entire function, then it is known
that generically there exists a sequence $(p_{n},q_{n})_{n=0,1,2,\ldots}\;\; q_{n}>0$
, so that $lim_{n\rightarrow\infty}\Bigl[p_{n}/q_{n}\Bigr]_{f}=\mathcal{\mathcal{\mathcal{\mathtt{\mathrm{\,\mathfrak{f}}}}}}$
(\cite{r8}).

In the present paper we investigate all the possible limits of sequences
$\Bigl(\Bigl[p_{n}/q_{n}\Bigr]_{f}\Bigr)_{n=0,1,2,..}\;$ on compact
subsets $K$ of $\mathbb{C}\setminus\{0\}$ or compact sets disjoint
from the domain of definition of the holomorphic function $\mathcal{\mathtt{\mathrm{\mathcal{\mathfrak{f}}}}}$.
We show that generically all functions holomorphic in a neighborhood
of $K$ are such limits. Thus, we have formal power series (or holomorphic
functions on a domain $\Omega$) with universal Pad\'{e} approximants.
The particular case $q=0$ is the well known case of universal Taylor
series (\cite{r6}, \cite{r16}, \cite{r32}, \cite{r33}, \cite{r35}, \cite{r37}, \cite{r41}, \cite{r42} and \cite{r49}) where the approximation
is realized by the partial sums. However, now we impose several conditions
on the approximating integers $(p_{n},q_{n})$ and $q_{n}$ may be
different from $0$; in particular, we can have universal approximations
with $p_{n}=q_{n}$, or $lim_{n\rightarrow\infty}p_{n}=lim_{n\rightarrow\infty}q_{n}=+\infty$
or $lim_{n\rightarrow\infty}(p_{n}-q_{n})=+\infty$ and others.

It was during the inspiring talk of Professor Brezinski that we got
the idea that the results on universal Taylor series may be transferred
to the case of Pad\'{e} approximants. Thus, we obtain the universal approximation
by Pad\'{e} approximants on compact sets $K$ with connected complement
$K^{c}$. However, the fact that Pad\'{e} approximants may also have poles
allows us to do approximation on compact sets $K$ of arbitrary connectivity
and these results are generic on spaces of holomorphic functions defined
on arbitrary planar domains and not just on simply connected domains.
This is not possible in the case of universal Taylor series where
the approximation is realized by polynomials (the partial sums). As
methods of proofs we use Baire\textquoteright{}s Category Theorem
combined with Runge\textquoteright{}s or Mergelyan\textquoteright{}s
Theorems. For the role of Baire\textquoteright{}s Category Theorem
in analysis we refer to \cite{r28} and \cite{r31} .

\section{PRELIMINARIES}

If $f={\displaystyle \sum_{v=0}^{\infty}a_{v}z^{v}}\in\mathbb{C^{\mathbb{N_{\textrm{0}}}}}$
and $p,\: q$ are non negative integers, we denote by\[
\Bigl[p/q\Bigr]_{f}\left(z\right)\]
a rational function of the form \[
\frac{{\textstyle \sum_{v=0}^{p}}n_{v}z^{v}}{\sum_{v=0}^{q}d_{v}z^{n}},\; d_{0}=1,\; n_{p}d_{q}\neq0\]
such that its Taylor development ${\textstyle {\displaystyle \sum_{v=0}^{\infty}}b_{v}z^{v}}$
coincides with ${\textstyle {\displaystyle \sum_{v=0}^{\infty}}a_{v}z^{v}}$up
to the first $p+q+1$ terms:\[
b_{v}=a_{v\;}forall\; v=0,1,\ldots,p+q.\]

Such a rational function does not always exist (see, for instance
\cite{r2} and \cite{r14}). It may also happen that there exist several such
rational functions (\cite{r22} and \cite{r26}). When such a rational fraction
exists, it is called a \textbf{Pad\'{e} approximant} of type $\left(p,q\right)$
to the series $f$. A necessary and sufficient condition for existence
and uniqueness is that the determinant of the \textbf{Hankel matrix}
$H^{\left(f\right)} _{q}\left(a_{p-q+1}\right)$of order $q$ at
$a_{p-q+1}$is different from zero:\[
\det\left(H_{q}^{\left(f\right)}\left(a_{p-q+1}\right)\right):=\det\left(\underbrace{\begin{array}{ccccc}
a_{p-q+1} & a_{p-q+2} & a_{p-q+3} & \ldots & a_{p}\\
a_{p-q+2} & a_{p-q+3} & a_{p-q+4} & \ldots & a_{p+1}\\
a_{p-q+3} & a_{p-q+4} & a_{p-q+5} & \ldots & a_{p+2}\\
\vdots & \vdots & \vdots & \vdots & \vdots\\
a_{p} & a_{p+1} & a_{p+2} & \ldots & a_{p+q-1}\end{array}}_{q}\right)\neq0.\]

Then we write $f\mathfrak{\in D}_{p,q}$ . If $f\mathfrak{\in D}_{p,q}$,
then the \textbf{Jacobi explicit formula} for $\Bigl[p/q\Bigr]_{f}\left(z\right)$
involves polynomial expressions on a finite number of the coefficients
$a_{v}$ of $f$ and its partial sums (\cite{r14}):\[
\left[p/q\right]_{f}\left(z\right)=\frac{\det\left(\begin{array}{cccc}
z^{q}S_{p-q}\left(z\right) & z^{q-1}S_{p-q+1}\left(z\right) & \cdots & S_{p}\left(z\right)\\
a_{p-q+1} & a_{p-q+2} & \cdots & a_{p+1}\\
\vdots & \vdots & \vdots & \vdots\\
a_{p} & a_{p+1} & \cdots & a_{p+q}\end{array}\right)}{\det\left(\begin{array}{cccc}
z^{q} & z^{q-1} & \cdots & 1\\
a_{p-q+1} & a_{p-q+2} & \cdots & a_{p+1}\\
\vdots & \vdots & \vdots & \vdots\\
a_{p} & a_{p+1} & \cdots & a_{p+q}\end{array}\right)}\]
with
\[ S_{k}\left(z\right) = \left\{ \begin{array}{ll}
		\sum_{v=0}^{k}a_{v}z^{v} & \mbox{for $k\geq0$,} \\
		0 & \mbox{for: $k<0$.}
		\end{array}
	\right.
\]

Thus, when we restrict our attention to $f$ \textquoteright{}s in
$\mathfrak{D}_{p,q},$ which is an open dense subset of the vector
space $C^{\mathbb{N_{\textrm{0}}}}$ endowed with the Cartesian topology,
then $\Bigl[p/q\Bigr]_{f}\left(z\right)$ varies continuously with
$f$.

We shall need the following lemmas.

\textbf{Lemma 2.1 }Let $F\left(z\right)$ and $G\left(z\right)$ be
two complex valued functions on a set $E$. We assume that $\sup_{z\in E}\mid F\left(z\right)\mid<+\infty,$$\inf_{z\in E}\mid G\left(z\right)\mid>0$
and $\sup_{z\in E}\mid G\left(z\right)\mid<+\infty$. Let also $\Phi\left(z,x\right)$
and $W\left(z,x\right)$ be two complex valued functions defined on
on $E\times\left(\mathbb{C}^{n}\setminus\left\{ 0\right\} \right)$
such that $\Phi\left(z,x\right)\rightarrow F\left(z\right)$ and $W\left(z,x\right)\rightarrow G\left(z\right)$
uniformly on $E$ as $x\rightarrow0$. Then, for every $a>0$ there
exists a ${\displaystyle \widetilde{\delta}>0}$ so that \[
\sup_{z\in E}\mid\frac{\Phi\left(z,x\right)}{W\left(z,x\right)}-\frac{F\left(z\right)}{G\left(z\right)}\mid<a\]

whenever $x\in\mathbb{C}^{n},\quad0<\left\Vert x\right\Vert <\widetilde{\delta}$
.$\blacksquare$

The proof is elementary and is omitted.

\textbf{Lemma 2.2} Let $M=M\left(d\right)=\left(a_{i,j}\left(d\right)\right)_{i,j=1}^{q}$
be a quadratic matrix with entries $a_{i,j}=a_{i,j}\left(d\right)$
depending linearly on a parameter $d$, in the sense that $a_{i,j}\left(d\right)=c_{i,j}d+\tau_{i,j}$.
Assume that
\[
\left\{ \begin{array}{ll}
		i+j<q\Rightarrow c_{i,j}=0 \\
		i+j=q\Rightarrow c_{i,j}=c\neq0\; and\:\tau_{i,j}=\tau .
	\end{array}
\right.
\]

Then the determinant $\det\left(M\right)$ of $M$ is a polynomial
in $d$ of degree $q$ with leading coefficient $\left(-1\right)^{q+1}c^{q}\neq0$.$\blacksquare$

The proof is elementary and is omitted.

\textbf{Lemma 2.3} Let $P\left(z\right)={\displaystyle \sum_{v=0}^{N}\varepsilon_{v}z^{v}}$
be a non-zero polynomial, $K\subset\subset\mathbb{C}$ be a compact
set and $a>0$. Suppose $p$ and $q$ are positive integer numbers
such that $p>\deg P+q$. Then there exist $c,\: d\in\mathbb{C\setminus}\left\{ 0\right\} $
such that the rational function\[
R\left(z\right)=P\left(z\right)+\frac{d\: z^{p}}{1-\left(c\: z\right)^{q}}\]
satisfies the following

\textbf{i.} $1-\left(c\: z\right)^{q}\neq0$, for all $z\in K.$

\textbf{ii.} $\sup_{z\in K}\left|P\left(z\right)-R\left(z\right)\right|<a.$

\textbf{iii.} The Taylor expansion $R\left(z\right)={\displaystyle \sum_{v=0}^{\infty}\beta_{v}z^{v}}$
of $R\left(z\right)$ around $0$ satisfies\[
\beta_{v}=\varepsilon_{v}\;\: for\; all\; v\leq\deg P.\]

\textbf{iv.} $R\left(z\right)\in\mathfrak{D}_{p,q}$ and the Pad\'{e}
approximant $\left[p/q\right]_{R}\left(z\right)$ coincides with the
rational function $R\left(z\right)$.

\textbf{Proof }Let $M<+\infty$ be such that $\mid z\mid\leq M$ for
all $z\in K$. Then for $z\in K$ we have $\mid\left(cz\right)^{q}\mid\leq\mid c\mid^{q}M^{q}<1$,
provided $\mid c\mid\leq\left(1/M\right)$. This implies \textbf{i}.
Application of Lemma 2.1 for $z\in E=K\cup\left\{ 0\right\} ,\: x=\left(c,d\right)\in\mathbb{C}^{2},\: F\left(z\right)\equiv0,\: G\left(z\right)\equiv1,\:\Phi\left(z,x\right)\equiv dz^{p}$
and $W\left(z,x\right)\equiv1-\left(cz\right)^{q}$ shows that there
is $\widetilde{\delta}$ such that\[
{\displaystyle \sup_{z\in E}}\mid\frac{\Phi\left(z,x\right)}{W\left(z,x\right)}-\frac{F\left(z\right)}{G\left(z\right)}\mid={\displaystyle \sup_{z\in K}}\mid\frac{dz^{p}}{1-\left(cz\right)^{p}}\mid=\sup_{z\in K}\mid P\left(z\right)-R\left(z\right)\mid<a\]

whenever $\left\Vert \left(c,d\right)\right\Vert <\widetilde{\delta}$
. This proves \textbf{ii}.To obtain \textbf{iii}, notice that the
identity\[
\frac{1}{1-\left(cz\right)^{q}}=1+\left(cz\right)^{q}+\left(cz\right)^{2q}+\ldots\]
implies that the Taylor development ${\displaystyle \sum_{v=0}^{\infty}\beta_{v}z^{v}}$
of $R\left(z\right)$ around $0$ is $R\left(z\right)=P\left(z\right)+d\: z^{p}\left(1+\left(c\, z\right)^{p}+\left(c\, z\right)^{2q}+\ldots\right)$.
Since $p>\deg P,\: P\left(z\right)$ is a partial sum of $R$'s Taylor
expansion, and $\beta_{v}=\varepsilon_{v}$ for all $v\leq\deg P$.
Finally, to show \textbf{iv}, observe that the Taylor coefficient
$\beta_{v}$ does not depend on the parameter $d$ whenever $v<p$.
We also have $\beta_{p}=d$. Further, in the expression of $\beta_{v}$,
for $v>p$, we have $\beta_{v}=d\,\tau_{v}$where $\tau_{v}$ is independent
of $d$. By Lemma 2.2, the Hankel determinant $H_{q}^{\left(R\right)}\left(\beta_{p-q+1}\right)=0$
is a polynomial in $d$ of degree $q$ with leading coefficient $\left(-1\right)^{q+1}$.
Hence, equation $H_{q}^{\left(R\right)}\left(\beta_{p-q+1}\right)=0$
with unknown $d$ has a finite number of solutions. Choosing a $d\in\mathbb{C}$
so that $H_{q}^{\left(R\right)}\left(\beta_{p-q+1}\right)\neq0$,
we infer $R\left(z\right)\in\mathfrak{D}_{p,q}$. Since the constant
term in $R\left(z\right)$'s denominator equals $1$, the uniqueness
of the Pad\'{e} approximant of type $\left(p,q\right)$ to the Taylor
series $P\left(z\right)+d\: z^{p}\left(1+\left(c\, z\right)^{p}+\left(c\, z\right)^{2q}+\ldots\right)$
guarantees that the Pad\'{e} approximant $\left[p/q\right]_{R}\left(z\right)$
coincides with the rational function $R\left(z\right)$. The proof
is complete.$\blacksquare$

\textbf{Lemma 2.4} Let ${\displaystyle \widetilde{P}\left(z\right)={\displaystyle \sum_{v=0}^{N}\varepsilon_{v}z^{v}}}$
be a non-zero polynomial, $K\subset\subset\mathbb{C}$ be a compact
set and $a>0$. Suppose $p$ and $q$ are positive integer numbers
such that $p,\: q>\deg\widetilde{P}$. Then there exist $c,\: d\in\mathbb{C\setminus}\left\{ 0\right\} $
such that the rational function\[
R\left(z\right)=\frac{\widetilde{P}\left(z\right)+dz^{p}}{1-\left(cz\right)^{q}}\]
 satisfies the following:

\textbf{i.} $1-\left(cz\right)^{q}\neq0$ for all $z\in K$.

\textbf{ii.} $\sup_{z\in K}\mid\widetilde{P}\left(z\right)-R\left(z\right)\mid<a$.

\textbf{iii.} The Taylor expansion $R\left(z\right)={\displaystyle \sum_{v=0}^{\infty}\beta_{v}z^{v}}$
of $R\left(z\right)$ around $0$ satisfies\[
\beta_{v}=\varepsilon_{v\;}for\: all\: v\leq\deg\widetilde{P}.\]

\textbf{iv. }$R\left(z\right)\in\mathfrak{D}_{p,q}$ and the Pad\'{e}
approximant $\left[p/q\right]_{R}\left(z\right)$ coincides with the
rational function $R\left(z\right)$.

\textbf{Proof }Let $M<+\infty$ be such that $\mid z\mid\leq M$ for
all $z\in K$. Then for $z\in K$ we have $\mid\left(cz\right)^{q}\mid\leq\mid c\mid^{q}M^{q}<1$,
provided $\mid c\mid\leq\left(1/M\right)$. This implies \textbf{i}.
To prove \textbf{ii}, it is enough to apply Lemma 2.1 for $z\in E=K,\: x=\left(c,d\right)\in\mathbb{C}^{2}$
$F\left(z\right)\equiv0,\: G\left(z\right)\equiv1,\:\Phi\left(z,x\right)\equiv c^{q}\, z^{q\,}\widetilde{P}\left(z\right)+d\, z^{p}$
and $W\left(z,x\right)\equiv1-\left(cz\right)^{q}$. It follows that
there exists s $\widetilde{\delta}>0$  such that\[
\sup_{z\in E}\mid\frac{\Phi\left(z,x\right)}{W\left(z,x\right)}-\frac{F\left(z\right)}{G\left(z\right)}\mid=\sup_{z\in K}\mid\frac{c^{q}\, z^{q}\:\widetilde{P}\left(z\right)+d\, z^{p}}{1-\left(cz\right)^{q}}\mid=\sup_{z\in K}\mid\widetilde{P}\left(z\right)-R\left(z\right)\mid<a\]
whenever $\left\Vert \left(c,d\right)\right\Vert <\widetilde{\delta}$
and \textbf{ii} is proved. To obtain \textbf{iii}, notice that the
Taylor development ${\displaystyle \sum_{v=0}^{\infty}\beta_{v}z^{v}}$
of $R\left(z\right)$ around $0$ is\[
R\left(z\right)=\widetilde{P}\left(z\right)+d\, z^{p}+\left(c\, z\right)^{q}\:\widetilde{P}\left(z\right)+\left(c\, z\right)^{q}\, d\: z^{p}+\left(c\, z\right)^{2q}d\: z^{p}+\ldots.\]
Since $p,q>\deg\widetilde{P},\;\widetilde{P}\left(z\right)$ is a
partial sum of $R$'s Taylor expansion. Thus, $\beta_{v}=\epsilon_{v}$
for all $v\leq\deg\widetilde{P}$. Finally, to show \textbf{iv}, observe
that

\[
\beta_{v}=
	\left\{
		\begin{array}{ll}
			\lambda_{v}d+\tau_{v},\mbox{ with\: $\lambda_{v}$\: and\:$\tau_{v}$\: independent\: of\: $d$} \\
			\tau_{v}, \mbox{with\:$\tau_{v}$ \: independed\: of\: $ d$ .}
		\end{array}
	\right.
\]

For $v<p$ , it holds $\beta_{v}=\tau_{v}$ and therefore, the coefficient
$\beta_{v}$ is independent of $d$. Further, if $q>0$ and $v=p$,
then $\beta_{v}=d+\tau_{v}$; and if $v>p$ and $v<p+q-1$, we have
$\beta_{v}=\tau_{v}+\lambda_{v}d$ where $\tau_{v}$ is independent
of $d$. Thus, by Lemma 2.2, the Hankel determinant $H_{q}^{\left(R\right)}\left(\beta_{p-q+1}\right)$
is a polynomial in $d$ of degree $q$ with leading coefficient $\left(-1\right)^{q+1}$.
Hence, equation $H_{q}^{\left(R\right)}\left(\beta_{p-q+1}\right)=0$
with unknown $d$ has a finite number of solutions and we can choose
$d$ so that $H_{q}^{\left(R\right)}\left(\beta_{p-q+1}\right)\neq0$.
We infer $R\left(z\right)\in\mathfrak{D}_{p,q}$. Since the constant
term in $R\left(z\right)$'s denominator equals $1$, the uniqueness
of the Pad\'{e} approximant of type $\left(p,q\right)$ to the Taylor
series $\widetilde{P}\left(z\right)+d\, z^{p}+\left(c\, z\right)^{q}\,\widetilde{P}\left(z\right)+\left(c\, z\right)^{q}\, d\, z^{p}+\left(c\, z\right)^{2q}\widetilde{P}\left(z\right)+\left(c\, z\right)^{2q}\, d\, z^{p}+\ldots$
guarantees that the Pad\'{e} approximant $\left[p/q\right]_{R}\left(z\right)$
coincides with the rational function $R\left(z\right)$. The proof
is complete. $\blacksquare$

\textbf{Lemma 2.5} Let $A\left(z\right)$, $B\left(z\right)$ and
$P\left(z\right)={\displaystyle \sum_{\nu=0}^{N}\varepsilon_{\nu}z^{\nu}}$
be nonzero polynomials. Let also $K\subset\subset\mbox{\ensuremath{\mathbb{C}}}$
be a compact set, $\lambda$ be a positive integer such that $\lambda>\deg\mbox{\ensuremath{P}}$
and\textbf{ $a>0$.} Assume that $B\left(z\right)\neq0$ for all $z\in K$
and $B\left(0\right)=1$. Then, for any $q\geq\deg B$ and $p>\max\left\{ \lambda+degA,\: q+degP\right\} $,
there exist two non zero polynomials $\tilde{A}\left(z\right)$ and
$\tilde{B}\left(z\right)$, with degrees $p=\deg\tilde{A}$ and $q=\deg\tilde{B}$
respectively, such that the rational function

\[
R\left(z\right)=\frac{\tilde{A}\left(z\right)}{\tilde{B}\left(z\right)}\]
 satisfies the following.

\textbf{i.} $\tilde{B}\left(z\right)\neq0$, for all $z$$\in K\cup\left\{ 0\right\} .$

\textbf{ii.} $\sup_{z\in K}\left|R\left(z\right)-\left[P\left(z\right)+z^{\lambda}\frac{A\left(z\right)}{B\left(z\right)}\right]\right|$$<a$.

\textbf{iii.} The Taylor expansion $R\left(z\right)=$${\displaystyle \sum_{\nu=0}^{\infty}\beta_{\nu}z^{\nu}}$of
$R\left(z\right)$ around $0$ satisfies\[
\beta_{\nu}=\varepsilon_{\nu},\: for\; all\;\nu\leq\deg P.\]
 \textbf{iv. }$R\left(z\right)\in\mathfrak{D}_{p,q}$, and the Pad\'{e}
approximant $\left[p/q\right]_{R}\left(z\right)$ coincides with the
rational function $R\left(z\right)$.

\textbf{Proof }Let $M<+\infty$ be such that $\left|z\right|\leq M$
for all $z\in K.$ Put \[
\tilde{A}\left(z\right)\equiv\left[B\left(z\right)+c\, z^{q}\right]\; P\left(z\right)+z^{\lambda}\, A\left(z\right)+d\, z^{p},\;\tilde{B}\left(z\right)\equiv B\left(z\right)+cz^{q}\]
 and

\[
R\left(z\right)=\frac{\tilde{A}\left(z\right)}{\tilde{B}\left(z\right)}.\]

Then for $z\in K$ we have $\left|\tilde{B}\left(z\right)\right|\geq\inf{}_{z\in K\cup\left\{ 0\right\} }\left|B\left(z\right)\right|-\left|c\right|M^{q}>0$,
provided $\left|c\right|<\left(\inf{}_{z\in K\cup\left\{ 0\right\} }\left|B\left(z\right)\right|/M^{q}\right)$.
This implies \textbf{i}. Further, application of Lemma 2.1 for $z\in E=K,\; x=\left(c,d\right)\in\mathbb{C}^{2}$,
$F\left(z\right)\equiv0$, $G\left(z\right)\equiv B^{2}\left(z\right)$,
$\Phi\left(z,x\right)\equiv-cz^{q+\lambda}A\left(z\right)+dz^{p}B\left(z\right)$
and $W\left(z,x\right)\equiv\left[B\left(z\right)\right]\left[B\left(z\right)+cz^{q}\right]$
shows that there is a $\tilde{\delta}$ such that\[
\sup_{z\in E}\left|\frac{\Phi\left(z,x\right)}{W\left(z,x\right)}-\frac{F\left(z\right)}{G\left(z\right)}\right|=\sup_{z\in K}\left|\frac{z^{\lambda}A\left(z\right)}{B\left(z\right)+cz^{q}}-\frac{z^{\lambda}A\left(z\right)}{B\left(z\right)}+\frac{dz^{p}}{B\left(z\right)+cz^{q}}\right|=\]

\[
=\sup\left|R\left(z\right)-\left[P\left(z\right)+z^{\lambda}\frac{A\left(z\right)}{B\left(z\right)}\right]\right|<a,\]
whenever $\left\Vert \left(c,d\right)\right\Vert <\widetilde{\delta}$.
This proves \textbf{ii}. To prove \textbf{iii}, we note that

\[
\left[B\left(z\right)+cz^{q}\right]^{-1}=\left[B\left(0\right)\right]^{-1}+\tau_{1}z+\tau_{2}z^{2}+\ldots\]
and therefore the Taylor development ${\displaystyle \sum_{\nu=0}^{\infty}\beta_{\nu}z^{\nu}}$of
$R\left(z\right)$ around $0$ is

\[
R\left(z\right)=P\left(z\right)+\left(\frac{z^{\lambda}A\left(z\right)}{B\left(0\right)}+\tau_{1}z^{\lambda+1}A\left(z\right)+\ldots\right)+\left(\frac{dz^{p}}{B\left(0\right)}+\tau_{1}dz^{p+1}+\tau_{2}dz^{p+2}+\ldots\right).\]

Since $\lambda>\deg\mbox{\ensuremath{P}}$, the polynomial $P\left(z\right)={\displaystyle \sum_{\nu=0}^{N}\varepsilon_{\nu}z^{\nu}}$
is a partial sum of $R$'s Taylor expansion. Thus, $\beta_{\nu}=\varepsilon_{\nu}$
whenever $\nu\leq\deg P$. It remains to show \textbf{iv}. To do so,
remind that

$R\left(z\right)=P\left(z\right)+\left(\left[B\left(0\right)\right]^{-1}z^{\lambda}A\left(z\right)+\tau_{1}z^{\lambda+1}A\left(z\right)+\tau_{2}z^{\lambda+2}A\left(z\right)+\ldots\right)+$

$+\left(\left[B\left(0\right)\right]^{-1}dz^{p}+\tau_{1}dz^{p+1}+\tau_{2}dz^{p+2}+\ldots\right)$.

So, for any $\nu<p$, the parameter $d$ does not appear in $R$\textquoteright{}s
Taylor coefficient $\beta_{\nu}$. For $\nu=p$, it holds, $\beta_{p}=\left(\left[d/B\left(0\right)+w_{p}\right]\right)$,
with $w_{p}$ independent of $d$. If $\nu>p$, the Taylor coefficient
$\beta_{\nu}$ depends at most linearly on $d$. From Lemma 2.2, it
follows that the Hankel determinant $H_{q}^{\left(R\right)}\left(\beta_{p-q+1}\right)$
is a polynomial in $d$ of degree $q$ with leading coefficient $\left(-1\right)^{q+1}$.
This determinant vanishes on a finite set of values of $d$. We can
avoid these values and choose $d$ so that $H_{q}^{\left(R\right)}\left(\beta_{p-q+1}\right)\neq0$.
We infer $R\left(z\right)\in\mathfrak{D}_{p,q}$. Since the constant
term in $R\left(z\right)$\textquoteright{}s denominator equals $1$,
the uniqueness of the Pad\'{e} approximant of type $\left(p,q\right)$
to the above Taylor series guarantees that the Pad\'{e} approximant $\left[p/q\right]_{R}\left(z\right)$
coincides with the rational function $R\left(z\right)$. The proof
is complete.$\blacksquare$

\paragraph{Remark 2.6 \textmd{In the above lemmas the rational function $R\left(z\right)$
is always a quotient of a polynomial of degree $p$ as numerator and
of a polynomial of degree $q$ as denominator. With a little more
effort, we can guarantee that these polynomials don\textquoteright{}t
have common zeros in $\mathbb{C}$. }}

\paragraph{Lemma 2.7 \textmd{Let $K$ be a compact subset of $\mathbb{C}$ such
that $0\notin K$ and the complement $K^{c}$ of $K$ is connected.
Let $\bar{D}$ be a closed disk centered at $0$ with radius $r>0$
and such that $\overline{D}\bigcap K=\slashed{O}$. Let also $P\left(z\right)$
and $Q\left(z\right)$ be two analytic polynomials. Then, for any
$\varepsilon>0$, there exist a polynomial $\tilde{P}\left(z\right)$
such that }}

\textbf{i.} $\sup_{z\in\bar{D}}\left|\tilde{P}\left(z\right)-P\left(z\right)\right|<\varepsilon$.

\textbf{ii.} $\sup_{z\in K}\left|\tilde{P}\left(z\right)-Q\left(z\right)\right|<\varepsilon$.

\textbf{iii.} A partial sum of $\tilde{P}\left(z\right)$ is exactly
the polynomial $P\left(z\right)$.

\textbf{Proof }Let $\lambda>\deg P$. We are looking for $\tilde{P}\left(z\right)=P\left(z\right)+z^{\lambda}\Pi\left(z\right)$,
where $\Pi\left(z\right)$ is a polynomial suitably chosen. To do
so, it is enough to have\[
\sup_{z\in K}\left|\Pi\left(z\right)-\frac{Q\left(z\right)-P\left(z\right)}{z^{\lambda}}\right|<\frac{\varepsilon}{M^{\lambda}}\;\, and\;\,\sup_{z\in\bar{D}}\left|\Pi\left(z\right)\right|<\frac{\varepsilon}{r^{\lambda}},\]
 where $M=\sup_{z\in K}\left|z\right|$. Indeed, it suffices to approximate
on $K\cup\bar{D}$ the function $\mathfrak{f}$ defined by

\[
\mathfrak{f}(z)=\left\{ \begin{array}{cc}
0 & on\,\bar{D}\\
\frac{Q\left(z\right)-P\left(z\right)}{z^{\lambda}} & on\, K\end{array},\right.\]

This is possible by Runge's theorem. The proof is complete.$\blacksquare$

\section{SELEZNEV- PAD\'{E} UNIVERSAL APPROXIMANTS}

\textbf{Definition 3.1 }Let $\mathfrak{F}\subset\mathbb{N_{\textrm{0}}^{\textrm{2}}}$$\equiv\left\{ 0,1,2,\ldots\right\} \times\left\{ 0,1,2,\ldots\right\} $.

\textbf{(i).} We say that $\mathfrak{\mathfrak{F}}$ satisfies condition
$\left(\mathfrak{\mathcal{H}}\right)$, if it contains a sequence
$\left(p_{n},q_{n}\right)_{n=1,2,\ldots}\in\mathfrak{F}$ satisfying
at least one of the following conditions:

$\left(\mathcal{H_{\textrm{1}}}\right)\lim_{n\rightarrow\infty}p_{n}=\lim_{n\rightarrow+\infty}q_{n}=+\infty.$

$\left(\mathcal{H_{\textrm{2}}}\right):\lim_{n\rightarrow\infty}\left(p_{n}-q_{n}\right)=+\infty.$

\textbf{(ii).} We denote by \[
\mathcal{U}\]

the class of all formal power series $f={\displaystyle \sum_{\nu=0}^{\infty}a_{\nu}z^{\nu}}$satisfying
the following condition:

For every compact set $K\in\mathbb{C}\setminus\left\{ 0\right\} $
with connected complement $K^{c}$ and every function $h:K\rightarrow\mathbb{C}$
which is continuous on $K$ and holomorphic in the interior $intK$
of $K\;\left(:h\in A\left(K\right)\right)$, there exists a sequence
$\left(\tilde{p}_{n},\tilde{q}_{n}\right)_{n=1,2,\ldots}\in\mathfrak{F}$
such that

\textbf{(a).} $f\in\mathfrak{D}_{\tilde{p}_{n},\tilde{q}_{n}}$ for
all $n=1,2,\ldots$ and

\textbf{(b).} $\lim_{n\rightarrow\infty}\left[\tilde{p}_{n}/\tilde{q}_{n}\right]_{f}\left(z\right)=h\left(z\right)$
uniformly on $K$.

\paragraph{Remark 3.2 \textmd{In the above definition of $\mathfrak{\mathcal{U}}$,
condition (b) is equivalent to require $\lim_{n\rightarrow\infty}\left(d^{l}\left[\tilde{p}_{n}/\tilde{q}_{n}\right]_{f}/dz^{l}\right)\left(z\right)=h^{\left(l\right)}\left(z\right)$
uniformly on $K$, for any $l=0,1,2,...$ and $h$ every polynomial.
This is due to the fact that $\mathbb{C\setminus}\left\{ 0\right\} $
is open; see \cite{r17}. }}

\textbf{Remark 3.3 }If the sequence $\left(\tilde{p}_{n},\tilde{q}_{n}\right)_{n=1,2,\ldots}\in\mathfrak{F}$
takes infinitely many values, then we can pass to a subsequence which
takes every value at most once. If the sequence $\left(\tilde{p}_{n},\tilde{q}_{n}\right)_{n=1,2,\ldots}\in\mathfrak{F}$
takes finitely many values, then \[
h\left(z\right)=\left[p/q\right]_{f}\left(z\right)\; for\: some\;\left(p,q\right)\in\left\{ 0,1,2,\ldots\right\} \times\left\{ 0,1,2,\ldots\right\} .\]
The set of all these $h$\textquoteright{}s is denumerable; so, there
is a sequence $\left(\varepsilon_{n}\right)_{n=1,2,\ldots}$such that
$\lim_{n\rightarrow\infty}\varepsilon_{n}=0$ and $h+\varepsilon_{n}\neq\left[p/q\right]_{f}$
for every $n=1,2,\ldots$ and $\left(p,q\right)\in\mathbb{N_{\textrm{0}}^{\textrm{2}}}$.
Thus, if $\left(\tilde{\tilde{p}}_{1},\tilde{\tilde{q}}_{1}\right)$,$\left(\tilde{\tilde{p}}_{2},\tilde{\tilde{q}}_{2}\right)$,...,$\left(\tilde{\tilde{p}}_{n-1},\tilde{\tilde{q}}_{n-1}\right)$
are already defined, we can choose $\left(\tilde{\tilde{p}}_{n},\tilde{\tilde{q}}_{n}\right)\in\mathfrak{F}$
different from $\left(\tilde{\tilde{p}}_{1},\tilde{\tilde{q}}_{1}\right)$,$\left(\tilde{\tilde{p}}_{2},\tilde{\tilde{q}}_{2}\right)$,...,$\left(\tilde{\tilde{p}}_{n-1},\tilde{\tilde{q}}_{n-1}\right)$
so that

\[
\sup_{z\in K}\left|\left[\tilde{\tilde{p}}_{n}/\tilde{q}\right]_{f}\left(z\right)-\left(h\left(z\right)+\varepsilon_{n}\right)\right|<\frac{1}{n}.\]
It follows that in Definition 3.1.(ii) it is equivalent if we require
in addition \[
\left(\tilde{p}_{n},\tilde{q}_{n}\right)\neq\left(\tilde{p}_{m},\tilde{q}_{m}\right)\; forall\; n,\, m\; with\; n\neq m.\]

We are now in position to formulate the first main result of the section.

\paragraph{Theorem 3.4 \textmd{The class $\mathcal{U}$ is dense and $G_{\delta}$
in the space $\mathbb{C^{\mathbb{N_{\textrm{0}}}}}$ of all formal
power series endowed with the Cartesian topology, provided $\mathfrak{F}$
satisfies condition $\left(\mathcal{H}\right)$. }}

The proof of Theorem 3.4 requires the following well known result.

\paragraph{Lemma 3.5 \textmd{(\cite{r35}, \cite{r37} and \cite{r49}) There exists a sequence $\left(K_{m}\subset\subset\mathbb{C}\setminus\left\{ 0\right\} \right)_{m=1,2,\ldots}$
of compact sets with connected complement $K_{m}^{c}$ such that the
following holds: }}

For every compact subset $K$ of $\mathbb{C}\setminus\left\{ 0\right\} $with
connected complement $K^{c}$, there exists a $m\in\left\{ 1,2,\ldots\right\} $,
such that $K\subset K_{m}$.$\blacksquare$

\paragraph{Proof of Theorem 3.4 \textmd{Regarding Definition 3.1 it is equivalent
to consider approximations only on the compact sets $K_{m}$. Because
if $h\in A\left(K\right)$ and $\varepsilon>0$ are given, then by
Mergelyan\textquoteright{}s Theorem, we can find a polynomial $Q\left(z\right)$
such that} }

\[
\sup_{z\in K}\left|Q\left(z\right)-h\left(z\right)\right|<\frac{\varepsilon}{2}.\]

Let $m\in\left\{ 1,2,\ldots\right\} $ be such that $K\subset K_{m}$
(see Lemma 3.5). Then, $Q\left(z\right)$ being a polynomial belongs
to $A\left(K_{m}\right)$ and can be approximated on $K_{m}$ by a
Pad\'{e} approximant $\left[p/q\right]_{R}\left(z\right)$, with $\left(p,q\right)\in\mathfrak{F}$:
\[
\sup_{z\in K_{m}}\left|\left[p/q\right]_{f}\left(z\right)-Q\left(z\right)\right|<\frac{\varepsilon}{2}.\]
Thus, the triangle inequality implies:

\[
\sup_{z\in K}\left|\left[p/q\right]_{f}\left(z\right)-h\left(z\right)\right|<\frac{\varepsilon}{2}+\frac{\varepsilon}{2}=\varepsilon.\]
Further, by Mergelyan\textquoteright{}s Theorem, it suffices to approximate
not all functions $h\in A\left(K_{m}\right)$ but only the polynomials,
with coefficients in $\mathbb{Q}+i\mathbb{Q}$; the set of these polynomials
is countable and can be represented as a sequence \[
\left(Q_{j}\right)_{j=1,2,\ldots}.\]
Next, we define the sets\[
\mathbb{E}\left(m,\, j,\, s,\,\left(p,q\right)\right)=\left\{ f\in\mathbb{C^{\mathbb{N_{\textrm{0}}}}}:\: f\in\mathfrak{D}_{p,q}\, and\,\sup_{z\in K_{m}}\left|\left[p/q\right]_{f}\left(z\right)-Q_{j}\left(z\right)\right|<\frac{1}{s}\right\} .\]

It is easily seen that

\[
\mathcal{U}=\bigcap_{m,j,s=1,2,\ldots}\bigcup_{\left(p,q\right)\in\mathfrak{F}}\mathbb{E}\left(m,j,s,\left(p,q\right)\right).\]
(See \cite{r37} and \cite{r41}).

To complete the proof of Theorem 3.4, we also need the following.

\textbf{Lemma 3.6} For any $m,\, j,\, s\in\left\{ 1,2,\ldots\right\} $
and any $\left(p,q\right)\in\left\{ 0,1,2,\ldots\right\} \times\left\{ 0,1,2,\ldots\right\} $,
the set $\mathbb{E}\left(m,j,s,\left(p,q\right)\right)$is open in
$\mathbb{C}^{\mathbb{N_{\textrm{0}}}}$.

\textbf{Proof }This follows from the fact that $\mathfrak{D}_{p,q}$
is open in $\mathbb{C^{\mathbb{N_{\textrm{0}}}}}$ and that, by Jacobi
explicit formula, the coefficients of the numerator and the denominator
of $\left[p/q\right]_{f}\left(z\right)$ $\left(d_{0}=1\right)$vary
continuously with $f\in\mathfrak{D}_{p,q}$.$\blacksquare$

\textbf{Proposition 3.7 }For every $m,\, j,\, s\in\left\{ 1,2,\ldots\right\} $,
the set $\bigcup_{\left(p,q\right)\in\mathfrak{F}}\mathbb{E}\left(m,j,s,\left(p,q\right)\right)$
is open and dense in $\mathbb{C}^{\mathbb{N_{\textrm{0}}}}$.

\textbf{Proof }Lemma 3.6 implies that the set is open as a union of
open sets. In order to prove density, let $P$ be any non-zero polynomial.
It suffices to find $\left(p,q\right)\in\mathfrak{F}$ and $f\in\mathbb{E}\left(m,j,s,\left(p,q\right)\right)$
such that $P\left(z\right)$ is a partial sum of $f\left(z\right)$.
To do so, let us recall that $\mathfrak{F}$ satisfies condition $\left(\mathcal{H}\right)$;
that is $\mathfrak{F}$ contains a sequence $\left(p_{n},q_{n}\right)_{n=1,2,\ldots}\in\mathbb{N}_{0}^{2}$
satisfying condition $\left(\mathcal{H}_{1}\right)$ or condition
$\left(\mathcal{H}_{2}\right)$. Assume that $\left(\mathcal{H}_{1}\right)$
holds. Since $\lim_{n\rightarrow+\infty}p_{n}=\lim_{n\rightarrow+\infty}q_{n}=+\infty$,
there exists a $n_{0}$ so that $p_{n_{0}},\: q_{n_{0}}>degP$. We
set\[
p=p_{n_{0}},\: q=q_{n_{0}},\: K=K_{m}\:\; and\;\: a=\frac{1}{s}\]
and we apply Lemma 2.4. To do this, we consider a small closed disk
$\bar{D}$ centered at zero disjoint from $K_{m}$. We also consider
the function $w$ defined by $w\left(z\right)=P\left(z\right)$ on
$\bar{D}$ and $w\left(z\right)=Q_{j}\left(z\right)$ on $K_{m}$.
By Lemma 2.7 we find a polynomial $\tilde{P}\left(z\right)$ approximating
$w\left(z\right)$ on $K_{m}\cup\bar{D}$ and such that a partial
sum of $\tilde{P}\left(z\right)$ is $P\left(z\right)$. Then we find
a rational function $R\left(z\right)$ given by Lemma 2.4. Letting
$f=R\left(z\right)$, we infer that $f\in\mathbb{E}\left(m,\, j,\, s,\,\left(p,q\right)\right)$,
$\left(p,q\right)\in\mathfrak{F}$ and a partial sum of $f\left(z\right)$
is $P\left(z\right)$ as required. Assume now that $\left(\mathcal{H}_{2}\right)$
holds. Since $\lim_{n\rightarrow+\infty}\left(p_{n}-q_{n}\right)=+\infty$,
$q_{n}\geq0$, we find a $n_{0}$ so that $p_{n_{0}}-q_{n_{0}}>\deg P$.
Therefore, we apply Lemma 2.3 and we have the result as in the previous
case. This completes the proof of Proposition 3.7.$\blacksquare$

\textbf{End of Proof of Theorem 3.4} Since $\mathbb{C^{\mathbb{N_{\textrm{0}}}}}$
is a complete metric space, Baire\textquoteright{}s Category Theorem
combined with Proposition 3.7 implies that the class $\mathcal{U}$
is dense and $G_{\delta}$ in $\mathbb{C^{\mathbb{N_{\textrm{0}}}}}$.
This completes the proof of Theorem 3.4.$\blacksquare$

\textbf{Remark 3.8 }It is an open question if $\mathcal{U}$ contains
a dense vector space or closed infinite subspace except $0$ in $\mathbb{C^{\mathbb{N_{\textrm{0}}}}}$.

\textbf{Remark 3.9 }A careful examination of the proof of Theorem
3.4 shows that it remains also valid if we replace the Cartesian topology
with the topology of the ring $\mathbb{C}\left[\left[\mathbb{Z}\right]\right]$
of formal power series over $\mathbb{C}$. Recall that the distance
between two distinct sequences $\left(a_{\nu}\right)_{\nu\in\mathbb{N}}\in\mathbb{C}\left[\left[\mathbb{Z}\right]\right]$
and $\left(b_{\nu}\right)_{\nu\in\mathbb{N}}\in\mathbb{C}\left[\left[\mathbb{Z}\right]\right]$
in this topology is defined to be $dist\left(\left(a_{\nu}\right),\left(b_{\nu}\right)\right)=2^{-\kappa}$,
where $\kappa$ is the smallest natural number such that $a_{\kappa}\neq b_{\kappa}$.
The distance between two equal sequences is, of course, zero.

\textbf{Remark 3.10 }In the case $\mathfrak{F}=\left\{ \left(p,0\right):p=0,1,2,\ldots\right\} $
the result of Theorem 3.4 is a well known result of Seleznev (see
\cite{r6} and \cite{r49}).

\textbf{Remark 3.11 }The case $\left(p_{n}\right)_{n=1,2,\ldots}$
is bounded is not covered by Theorem 3.4. However for $p_{n}=0$,
$q_{n}=n\,\,\left(n=1,2,\ldots\right)$, the function $Q\left(z\right)=z-5$
can not be uniformly approximated on the compact set $K=\left\{ z:\left|z-5\right|\leq2\right\} $
by functions of the form $1/P_{n}\left(z\right)\,\,\left(n=1,2,\ldots\right)$
where $P_{n}$ are analytic polynomials. For, otherwise we would have
\[
\frac{1}{2}\leq\left|\frac{1}{P_{n}\left(z\right)}\right|\leq\frac{3}{2}\; for\; n\geq n_{0}\; on\;\left\{ z:\left|z-5\right|=1\right\} .\]
This implies $\left|P_{n}\left(z\right)\right|\leq2$ for all $n\geq n_{0}$
on $\left\{ z:\left|z-5\right|=1\right\} $. The maximum principle
implies $\left|P_{n}\left(5\right)\right|\leq2$ for all $n\geq n_{0}$.
Thus,\[
\left|\frac{1}{P_{n}\left(5\right)}\right|\geq\frac{1}{2},\; for\: all\; n\geq n_{0},\]
 which contradicts $\lim_{n\rightarrow\infty}\frac{1}{P_{n}\left(5\right)}=0$.

\textbf{Remark 3.12 }Cases covered by Theorem 3.4 are the following:

\textbf{(i).} $l>0$ and $p_{n}\geq lq_{n},\:\: q_{n}=n$,\textbf{ }

\textbf{(ii).} $l\in\mathbb{Q},\:\: l>0$, and $p_{n}=lq_{n}$, with
$q_{n}\rightarrow+\infty$.

\textbf{(iii).} $q_{n}\leq M$ and $p_{n}\rightarrow+\infty$. (The
case $M=0$ corresponds to Seleznev\textquoteright{}s result).

Especially, by \textbf{(ii)} of Remark 3.12 we see that we can also
work on the diagonal $p=q$ of the Pad\'{e} table.

\section{Approximation on arbitrary compact sets in $\mathbb{C}\backslash\{0\}$ }

Lemma 2.5 allows approximating holomorphic functions in a neighborhood
of a compact set $K\subset\mathbb{C}$ with arbitrary connectivity.

\textbf{Definition 4.1 }A set $\Im\subset\mathbb{N}_{0}^{2}\equiv\left\{ 0,1,2,...\right\} \times\left\{ 0,1,2,...\right\} $
satisfies condition $\left(\widetilde{\mathcal{H}}\right)$if it contains
a sequence $\left(p_{n},q_{n}\right)_{n=1,2,...}\in\Im$ such that:
\[
lim_{n\rightarrow\infty}q_{n}=lim_{n\rightarrow\infty}(p_{n}-q_{n})=+\infty.\blacksquare\]

Let $\left(K_{m}\right)_{m=1,2,...}$be a fixed sequence of compact
subsets of $\mathbb{C}\setminus\{0\}$.

\textbf{Definition 4.2 }Let $\Im\subset\mathbb{N}_{0}^{2}\equiv\left\{ 0,1,2,...\right\} \times\left\{ 0,1,2,...\right\} $
and $\left(K_{m}\right)_{m=1,2,...}$be as above. We denote by\[
\widetilde{\mathcal{U}}\]
the class of all formal power series $f={\displaystyle \sum_{v=0}^{\infty}a_{v}z^{v}\in\mathbb{C}^{\mathbb{N}{}_{0}}}$such
that for every $m=1,2,...$ and every function $Q\left(z\right)$
holomorphic in a neighborhood of $K_{m}$ the following holds.

There exists a sequence $\left(\widetilde{p_{n}},\widetilde{q_{n}}\right)_{n=1,2,\ldots}$
in $\Im$ such that

\textbf{(i).} $f{\displaystyle \in\mathfrak{D}_{\widetilde{p_{n}},\widetilde{q_{n}}}}$for
all $n=1,2,...$ and

\textbf{(ii).} $lim_{n\rightarrow\infty}\left[\widetilde{p_{n}}/\widetilde{q_{n}}\right]_{f}\left(z\right)=Q\left(z\right)$uniformly
on $K_{m}$.$\blacksquare$

It is easy to see that $\widetilde{\mathcal{U}}$ remains unchanged
if we require in addition \[
\left(\widetilde{p_{n}},\widetilde{q_{n}}\right)\neq\left(\widetilde{p_{l}},\widetilde{q_{l}}\right)for\; n\neq l.\]

\textbf{Remark 4.3 }Since $\mathbb{C}\setminus\{0\}$ is open, according
to \cite{r17} condition \textbf{(ii)} in the above Definition 4.2 is
equivalent to the requirement $lim_{n\rightarrow\infty}\left(d\,^{l}\:\left[\widetilde{p_{n}}/\widetilde{q_{n}}\right]_{f}/dz^{l}\right)\left(z\right)=Q^{\left(l\right)}\left(z\right)$uniformly
on $K_{m}$, for any $l=0,1,2,...$.

\textbf{Proposition 4.4 }If $\Im$ satisfies condition $(\widetilde{\mathcal{H}})$
and $\left(K_{m}\right)_{m=1,2,...}$are as above, then the class
$\widetilde{\mathcal{U}}$ is dense and $G_{\delta}$ in the space
$\mathbb{C^{\mathbb{N}}}^{_{_{0}}}$ of all formal power series endowed
with the Cartesian topology.

\textbf{Proof }By Runge\textquoteright{}s Theorem, it suffices to
approximate on each $K_{m}$ all rational functions with poles off
$K_{m}\cup\left\{ 0\right\} $. By choosing all coefficients of the
numerator and the denominator from $\mathbb{Q}+i\mathbb{Q}$, it is
easy to see that it is enough to approximate on $K_{m}$ a denumerable
set of rational functions $Q_{j,m}\;\left(j=1,2,...\right)$with poles
off $K_{m}\cup\left\{ 0\right\} $. Fix the functions\[
Q_{j,m}\;\left(j=1,2,...\: and\: m=1,2,...\right)\]
and consider the sets

\[
\mathbb{E}\left(m,j,s,\left(p,q\right)\right)=
	\left\{
		\begin{array}{c}
			f\in\mathbb{C^{\mathbb{N}}}^{_{_{0}}}:\; f\in\mathfrak{D}_{p,q}\: and \\
			\sup_{z\in K_{m}}\left|\left[p/q\right]_{f}\left(z\right)-Q_{j,m}\left(z\right)\right|<\frac{1}{s}
		\end{array}
	\right\}
\]

It is easily seen that\[
\widetilde{\mathcal{U}}={\displaystyle \bigcap_{m,j,s=1,2,\cdots}{\displaystyle \bigcup_{\left(p,q\right)\in\Im}\mathbb{E}\left(m,j,s,\left(p,q\right)\right)}}.\]

Next, each set $\mathbb{E}\left(m,\, j,\, s,\,\left(p,q\right)\right)$
is open. The justification is the same as that of Lemma 3.6. In order
to apply Baire\textquoteright{}s Category Theorem it remains to show
${\displaystyle \bigcup_{\left(p,q\right)\in\Im}\mathbb{E}\left(m,j,s,\left(p,q\right)\right)}$
is dense in $\mathbb{C^{N}}^{_{_{0}}}$. This can be done in a similar
way with the proof of Proposition 3.7. The only difference is that
we use Lemma 2.5 instead of Lemmas 2.3 and 2.4. In order to do this,
it is enough to approximate $Q_{j,m}\left(z\right)$ on $K_{m}$ by
a function of the form\[
P\left(z\right)+z^{\lambda}\frac{A\left(z\right)}{B\left(z\right)}\]
 with
\begin{itemize}
\item $P$ an arbitrary polynomial,
\item $\lambda>\deg P$,
\item $A,\: B$ arbitrary polynomials and
\item $B\left(z\right)\neq0$ on $K_{m}\cup\left\{ 0\right\} $.
\end{itemize}
For this purpose, it suffices that $A\left(z\right)/B\left(z\right)$
approximates $\left[Q_{j,m}\left(z\right)-P\left(z\right)\right]/z^{\lambda}$
on $K_{m}$. This is assured by Runge\textquoteright{}s Theorem. Since
we easily obtain $B\left(z\right)\neq0$ on $K_{m}\cup\left\{ 0\right\} $
, the result follows. $\blacksquare$

\textbf{Remark 4.5 }Cases covered by Theorem 4.4 are the following:

\textbf{(i).} $p_{n}\geq lq_{n}$, for all $n,\: l>1$ and $q_{n}\rightarrow+\infty$.

\textbf{(ii).} $l\in\mathbb{Q},\: l>1$ and $p_{n}=lq_{n}$, with
$q_{n}\rightarrow+\infty.$

A naturally posed question is can we have $p_{n}=q_{n}$?

We consider a non-empty finite union of open disks with centres in
$\mathbb{Q}+i\mathbb{Q}$ and rational radii and such that one of
these disks contains $0$. The set of all such unions is denumerable.
We denote by $\left(K_{m}\right){}_{m=1,2,\ldots}$the sequence of
complements of these unions in $\left\{ z\in\mathbb{C}:\left|z\right|\leqq n\right\} $,
where $n$ varies in the set of natural numbers. Then if a formal
power series $f={\displaystyle \sum_{v=0}^{\infty}a_{v}z^{v}\in\mathbb{C}^{\mathbb{N}{}_{0}}}$
belongs to the class $\widetilde{\mathcal{U}}$ (with respect to the
sequence $\left(K_{m}\right){}_{m=1,2,\ldots}$), it also has the
following property, provided $\Im$ satisfies condition $\widetilde{\mathcal{\left(H\right)}}$:
For every compact set $K\subset\mathbb{C}\setminus\left\{ 0\right\} $and
every function $Q$ holomorphic in a neighborhood of $K$ , the following
holds. There exists a sequence $\left(\widetilde{p_{n}},\widetilde{q_{n}}\right)_{n=1,2,\ldots}$
in $\Im$ such that

\textbf{(i).} $f{\displaystyle \in\mathfrak{D}_{\widetilde{p_{n}},\widetilde{q_{n}}}}$for
all $n=1,2,...$ and

\textbf{(ii). }$lim_{n\rightarrow\infty}\left[\widetilde{p_{n}}/\widetilde{q_{n}}\right]_{f}\left(z\right)=Q\left(z\right)$
uniformly on $K$.

The reason is that by Runge\textquoteright{}s Theorem we can approximate
$Q\left(z\right)$on $K$ by a rational function $\widetilde{Q}\left(z\right)$
with poles $r_{1},r_{2},\ldots,r_{N}$ in $\left[\left(\mathbb{Q}+i\mathbb{Q}\right)\setminus K\right]$.
Then we find $m=1,2,\ldots$ so that $K\subset K_{m}$ and $r_{1},r_{2},\ldots,r_{N}\in K_{m}^{c}$.
Thus, $\widetilde{Q}$ is holomorphic in a neighborhood of $K_{m}$
and can be approximated on $K_{m}$ by $\left[p/q\right]_{f}\left(z\right)$
as in Definition 4.2. Thus, we have proved the following.

\textbf{Theorem 4.6 }Let $\Im\subset\mathbb{N}_{0}^{2}\equiv\left\{ 0,1,2,\ldots\right\} \times\left\{ 0,1,2,\ldots\right\} $
satisfies condition $\widetilde{\mathcal{H}}$. Then there exists
a formal power series $f\in\mathbb{C}^{\mathbb{N}_{0}}$such that
the following holds. For every compact set $K\subset\mathbb{C}\setminus\left\{ 0\right\} $and
every function $Q\left(z\right)$ holomorphic in a neighborhood of
$K$, there exists a sequence $\left(\widetilde{p_{n}},\widetilde{q_{n}}\right)_{n=1,2,\ldots}\in\Im$
so that

\textbf{(i).} $f{\displaystyle \in\mathfrak{D}_{\widetilde{p_{n}},\widetilde{q_{n}}}}$for
all $n=1,2,...$ and

\textbf{(ii).} $lim_{n\rightarrow\infty}\left[\widetilde{p_{n}}/\widetilde{q_{n}}\right]_{f}\left(z\right)=Q\left(z\right)$uniformly
on $K$.

The collection of all such formal series $f$ is a dense $G_{\delta}$
subset of $\mathbb{C}^{\mathbb{N}_{0}}$endowed with the Cartesian
topology. It is also dense and $G_{\delta}$ in $\mathbb{C}\Bigl[\bigl[Z\bigr]\Bigr]$
endowed with the ring topology.$\blacksquare$

\section{SIMPLY CONNECTED DOMAINS}

Consider a simply connected domain $\Omega\neq\mathbb{C}$ containing
$0$. Let $f$ be a holomorphic function in $\Omega\left(:f\in\mathcal{O}\left(\Omega\right)\right)$
and let the Taylor development ${\displaystyle \sum_{v=0}^{\infty}}a_{v}z^{v}$
of $f$ around $0$. Let finally $p$ and $q$ be two non negative
integers. As it is already pointed out in §2, $f{\displaystyle \in\mathfrak{D}_{p,q}}$
if and only if $H_{q}^{\left(f\right)}\left(a_{p-q+1}\right)\neq0$.
Cauchy\textquoteright{}s estimates imply that $\mathfrak{D}_{p,q}$
is an open subset of $\mathcal{O}\left(\Omega\right)$ if $\mathcal{O}\left(\Omega\right)$
is endowed with the topology of uniform convergence on compact subsets
of $\Omega$.

Assume that $\left(L_{k}\right)_{k=1,2,\ldots}$is an increasing sequence
of compact subsets of $\overline{\Omega}$ such that
\begin{itemize}
\item $0\in int\left(L_{1}\right)$,
\item $\overline{L_{k}\cap\Omega}=L_{k}$ whenever $k=1,2,\ldots$,
\item $L_{k}^{c}$ is connected for all $k=1,2,\ldots$ and
\item every compact set $\mathcal{T}\subset\Omega$ is contained in some
$L_{k}$.
\end{itemize}
Let us consider the space $\mathcal{O}\left(\Omega,\left\{ L_{k}\right\} \right)$
of holomorphic functions $f\in\mathcal{O}\left(\Omega\right)$such
that, for each derivative $f^{\left(l\right)}\left(l=1,2,\ldots\right)$
and each $L_{k}\left(k=1,2,\ldots\right)$ the restriction\[
\mathit{f^{\left(l\right)}\mid L_{k}\cap\Omega}\]
is uniformly continuous on $L_{k}\cap\Omega$. Therefore, it extends
continuously on $\overline{L_{k}\cap\Omega}=L_{k}$. The space $\mathcal{O}\left(\Omega,\left\{ L_{k}\right\} \right)$
endowed with the seminorms\[
\sup_{z\in L_{k}}\left|f^{\left(l\right)}\left(z\right)\right|\;\left(l=1,2,\ldots\; and\; k=1,2,\ldots\right)\]
 becomes a Fréchet space containing the polynomials. Since we do not
know if the polynomials are dense in $\mathcal{O}\left(\Omega,\left\{ L_{k}\right\} \right)$,
we consider the closure $\overline{\mathbb{P}}\left(\Omega,\left\{ L_{k}\right\} \right)$
of the set of polynomials in $\mathcal{O}\left(\Omega,\left\{ L_{k}\right\} \right)$.
Then $\mathcal{O}\left(\Omega,\left\{ L_{k}\right\} \right)$ and
$\overline{\mathbb{P}}\left(\Omega,\left\{ L_{k}\right\} \right)$are
again Fréchet spaces and Baire\textquoteright{}s Category Theorem
is at our disposal. Let now $\left(K_{m}\right)_{m=1,2,\ldots}$be
a sequence of compact sets with
\begin{itemize}
\item connected complement $K_{m}^{c}$ and
\item $K_{m}\cap L_{k}=\textrm{Ø}$ whenever $\left(m=1,2,\ldots\; and\; k=1,2,\ldots\right).$
\end{itemize}
We are in position to prove the following:

\textbf{Theorem 5.1 }Suppose $\Im\subset\left\{ 0,1,2,\ldots\right\} \times\left\{ 0,1,2,\ldots\right\} $satisfies
condition $\left(\mathcal{H}\right)$. Under the above assumptions,
there exists a function $f\in\overline{\mathbb{P}}\left(\Omega,\left\{ L_{k}\right\} \right)$
such that for every polynomial $Q\left(z\right)$and every $m=1,2,\ldots$
there is a sequence $\left(p_{n},q_{n}\right)_{n=1,2,\ldots}\in\Im$
so that $f{\displaystyle \in\mathfrak{D}_{p_{n},q_{n}}\cap\overline{\mathbb{P}}\left(\Omega,\left\{ L_{k}\right\} \right)}$for
all $n$ and the following hold.

For every $l=0,1,2,\ldots$ we have
\begin{itemize}
\item $\sup_{z\in L_{k}}\left|\frac{d^{l}}{dz^{l}}\left[p_{n}/q_{n}\right]_{f}\left(z\right)-\frac{d^{l}}{dz^{l}}f\left(z\right)\right|\rightarrow0$
as $n\rightarrow\infty$ for all $k=1,2,\ldots$ and
\item $\sup_{z\in K_{m}}\left|\frac{d^{l}}{dz^{l}}\left[p_{n}/q_{n}\right]_{f}\left(z\right)-\frac{d^{l}}{dz^{l}}Q\left(z\right)\right|\rightarrow0$
as $n\rightarrow\infty$.
\end{itemize}
The set of these functions $f$ is dense and $G_{\delta}$ in the
space $\overline{\mathbb{P}}\left(\Omega,\left\{ L_{k}\right\} \right)$.

\textbf{Proof }It is easy to see that it suffices to approximate all
polynomials $Q_{j}\left(z\right)\;\left(j=1,2,\ldots\right)$ with
coefficients in $\mathbb{Q}+i\mathbb{Q}$ and their derivatives on
the compact sets $K_{m}\;\left(m=1,2,\ldots\right)$. We consider
the sets\[
\hat{\mathbb{E}}_{\Omega}\left(m,j,s,k,\left(p,q\right)\right)=\left\{ \begin{array}{c}
f\in\overline{\mathbb{P}}\left(\Omega,\left\{ L_{k}\right\} \right)\cap\mathfrak{D}_{p,q}:\\
\sup_{z\in K_{m}}\left|\frac{d^{l}}{dz^{l}}\left[p/q\right]_{f}\left(z\right)-\frac{d^{l}}{dz^{l}}Q_{j}\left(z\right)\right|<\frac{1}{s}\; and\\
\sup_{z\in L_{k}}\left|\frac{d^{l}}{dz^{l}}\left[p_{n}/q\right]_{f}\left(z\right)-\frac{d^{l}}{dz^{l}}f\left(z\right)\right|<\frac{1}{s}\\
for\; l=0,1,2,\ldots s\end{array}\right\} .\]

Obviously, the set of functions satisfying the hypotheses of the Theorem
can be written as\[
{\normalcolor \mathcal{\mathbf{\mathfrak{U}}}}\left(\Omega,\left\{ L_{k}\right\} ,\left\{ K_{m}\right\} \right):={\displaystyle \bigcap_{m,j,s,k=1}^{\infty}\bigcup_{\left(p,q\right)\in\Im}\hat{\mathbb{E}}_{\Omega}\left(m,\, j,\, s,\, k,\,\left(p,q\right)\right)}.\]

Cauchy estimates and the continuity of the operator ${\displaystyle \left(\mathfrak{D}_{p,q}\ni\right)f\longmapsto\left[p/q\right]_{f}}$
imply that $\hat{\mathbb{E}}_{\Omega}\left(m,\, j,\, s,\, k,\,\left(p,q\right)\right)$
is open in $\overline{\mathbb{P}}\left(\Omega,\left\{ L_{k}\right\} \right)$.
It remains to show that ${\displaystyle \bigcup_{\left(p,q\right)\in\Im}\hat{\mathbb{E}}_{\Omega}\left(m,\, j,\, s,\, k,\,\left(p,q\right)\right)}$
is dense in $\overline{\mathbb{P}}\left(\Omega,\left\{ L_{k}\right\} \right)$.
For this purpose fix any $L=L_{k}\;(k=1,2,\ldots)$ and consider a
polynomial $\varphi\left(z\right)$, $\varepsilon>0$ and $N\in\mathbb{N}\equiv\left\{ 1,2,\ldots\right\} $.
We are looking for a function $f\in\hat{\mathbb{E}}_{\Omega}\left(m,\, j,\, s,\, k,\,\left(p,q\right)\right)$
and a $\left(p,q\right)\in\Im$ such that
\begin{itemize}
\item $\sup_{z\in L}\left|\frac{d^{l}}{dz^{l}}f\left(z\right)-\frac{d^{l}}{dz^{l}}\varphi\left(z\right)\right|<\varepsilon$,
for $l=1,2,\ldots,N$
\item $\sup_{z\in K_{m}}\left|\frac{d^{l}}{dz^{l}}\left[p/q\right]_{f}\left(z\right)-\frac{d^{l}}{dz^{l}}Q_{j}\left(z\right)\right|<\frac{1}{s}$
and
\item $\sup_{z\in L}\left|\frac{d^{l}}{dz^{l}}\left[p/q\right]_{f}\left(z\right)-\frac{d^{l}}{dz^{l}}f\left(z\right)\right|<\frac{1}{s}$
for all $l=1,2,\ldots,s.$
\end{itemize}
By \cite{r17} and \cite{r27}, there exist two simply connected open sets\[
V\; with\; L\subset V\subset\Omega\; and\; W\; with\; K_{m}\subset W\]
such that $V\cap W=\textrm{Ø}.$ Considering an exhausting family
of compact sets for $W$, we find a compact set $S\subset W$ such
that $K_{m}\subset int\left(S\right)$ and the complement $S^{c}$
is connected. Similarly, we find a compact set $T$ such that $L\subset int\left(T\right)\subset T\subset V$.
We consider the function $w\left(z\right)=\varphi\left(z\right)$on
$V$ and the function $w\left(z\right)=Q_{j}\left(z\right)$ on $W$.
Runge\textquoteright{}s Theorem gives a polynomial $P_{n}\left(z\right)$
such that \[
\sup_{z\in T\bigcup S}\left|P_{n}\left(z\right)-w\left(z\right)\right|<\frac{1}{n}\;\left(n=1,2,\ldots\right).\]
Using Lemmas 2.3 and 2.4 we find rational functions $R_{n}\left(z\right)$
with poles outside $\Omega\cup S$ (close to $\infty$ ) such that
for some $\left(p,q\right)\in\Im$ depending on $n$
\begin{itemize}
\item $\sup_{z\in T\bigcup S}\left|R_{n}\left(z\right)-P_{n}\left(z\right)\right|<\frac{1}{n}$,
\item $R_{n}\left(z\right)\in\mathfrak{D}_{p,q}$ and
\item $\left[p/q\right]_{R_{n}}\left(z\right)\equiv R_{n}\left(z\right)$.
\end{itemize}
Since $lim_{n\rightarrow\infty}R_{n}\left(z\right)=Q_{j}\left(z\right)$
uniformly on the open set $int\left(S\right)\bigcup int\left(T\right)$,
Weierstrass Theorem implies that \[
lim_{n\rightarrow\infty}\left(d^{l}R_{n}\left(z\right)/dz^{l}\right)=\left(d^{l}w\left(z\right)/dz^{l}\right)\]
uniformly on $K_{m}\cup L\:\left(l=0,1,2,\ldots\right)$. Choosing
$f\left(z\right)$ to be one of the functions $R_{n}$$\left(z\right)$$\left(n\: big\right)$and
$\left(p,q\right)\in\Im$ as above, we find
\begin{itemize}
\item $\sup_{z\in L}\left|\frac{d^{l}}{dz^{l}}f\left(z\right)-\frac{d^{l}}{dz^{l}}\varphi\left(z\right)\right|<\varepsilon\;$$\left(l=1,2,\ldots,N\right)$,
\item $f\in\mathcal{O}\left(\Omega\right)\cap\mathfrak{D}_{p,q}$ and
\item $\sup_{z\in K_{m}}\left|\frac{d^{l}}{dz^{l}}f\left(z\right)-\frac{d^{l}}{dz^{l}}Q_{j}\left(z\right)\right|<\frac{1}{s}$
$\;$$\left(l=1,2,\ldots,s\right).$
\end{itemize}
Since $f\left(z\right)-\left[p/q\right]_{f}\left(z\right)\equiv0$
, Baire\textquoteright{}s Theorem completes the proof. $\blacksquare$

To give a first application of Theorem 5.1, let us consider
\begin{itemize}
\item a simply connected domain $\Omega$ containing $0$ and
\item an exhausting sequence $\left(L_{k}\right)_{k=1,2,\ldots}$ of compact
subsets of $\Omega$.
\end{itemize}
The sequence $\left(K_{m}\right)_{m=1,2,\ldots}$ is chosen to satisfy
the following.

\textbf{Lemma 5.2}(\cite{r37} and \cite{r40}) Let $\Omega$ be a domain of
$\mathbb{C}$. There exists a sequence $\left(K_{m}\right)_{m=1,2,\ldots}$
of compact subsets of $\mathbb{C}$ with $K_{m}\cap\Omega=\emptyset$
,$K_{m}^{c}$ connected, and such that the following holds.
\begin{itemize}
\item For every compact set $K\subset\mathbb{C}$, with $K\cap\Omega=\emptyset$
and $K^{c}$ connected, there exists a $m\in\left\{ 1,2,\ldots\right\} $so
that $K\subset K_{m}$ .$\blacksquare$
\end{itemize}
Then we obtain the following special case of Theorem 5.1.

\textbf{Theorem 5.3 }Let $\Omega$ be a simply connected domain of
$\mathbb{C}$ containing $0$. Let also $\Im\subset\left\{ 0,1,2,\ldots\right\} \times\left\{ 0,1,2,\ldots\right\} $satisfying
condition $\left(\mathcal{H}\right)$. There exists a holomorphic
function $f\left(z\right)\in\mathcal{O}\left(\Omega\right)$such that
\begin{itemize}
\item For every compact set $K\subset\mathbb{C}$, such that $K\cap\Omega=\emptyset$
and $K^{c}$ connected, and every polynomial $\varphi\left(z\right)$,
there exists a sequence $\left(p_{n},q_{n}\right)_{n=1,2,\ldots}\in\Im$
so that

\textbf{(i).} $f\in\mathfrak{D}_{p_{n},q_{n}}$ for all $n=1,2,\ldots$

\textbf{(ii).} $lim_{n\rightarrow\infty}\frac{d^{l}}{dz^{l}}\left[p_{n}/q_{n}\right]$$_{f}\left(z\right)=\frac{d^{l}}{dz^{l}}\varphi\left(z\right)$
uniformly on $K\quad\left(l=0,1,2,\ldots\right)$ and

\textbf{(iii).} $lim_{n\rightarrow\infty}\left[p_{n}/q_{n}\right]_{f}\left(z\right)=f\left(z\right)$,
uniformly on each compact subset of $\Omega$.

\textbf{(iv).} In particular, for every $h\left(z\right)\in A\left(K\right)$we
have $lim_{n\rightarrow\infty}\left[p_{n}/q_{n}\right]_{f}\left(z\right)=h\left(z\right)$,
uniformly on $K$.

\end{itemize}
The set of all such $f$ \textquoteright{} s is a dense $G_{\delta}$
subset of $\mathcal{O}\left(\Omega\right)$endowed with the topology
of uniform convergence on compacta. $\blacksquare$

To give a second application of Theorem 5.1, we consider an exhausting
sequence $\left(L_{k}\right)_{k=0,1,2,\ldots}$of compact subsets
of $\Omega$, with $0\in int\left(L_{0}\right)$. The sequence $\left(K_{m}\right)_{m=1,2,\ldots}$is
given by the following.

\textbf{Lemma 5.4}(\cite{r6}, \cite{r35}, \cite{r36} and \cite{r37}) Let $\Omega$ be a domain
of $\mathbb{C}$. There exists a sequence $\left(K_{m}\right)_{m=1,2,\ldots}$
of compact subsets of $\mathbb{C}$ with $K_{m}\bigcap\bar{\Omega}=\emptyset,$
$K_{m}^{c}$ connected, and such that the following holds.
\begin{itemize}
\item For every compact set $K\subset\mathbb{C},$ with $K\bigcap\bar{\Omega}=\mathbf{\mathrm{\slashed{O}}}$
and $K^{c}$ connected, there exists a $m\in\left\{ 1,2,...\right\} $so
that $K\subset K_{m}$. $\blacksquare$
\end{itemize}
Then we obtain the following special case of Theorem 5.1.

\textbf{Theorem 5.5 }Let $\Omega$ be a simply connected domain of
$\mathbb{C}$ containing $0$. Let also $\Im\subset\left\{ 0,1,2,\ldots\right\} \times\left\{ 0,1,2,\ldots\right\} $satisfying
condition $\left(\mathcal{H}\right)$. Then there exists a holomorphic
function $f\left(z\right)\in\mathcal{O}\left(\Omega\right)$such that
\begin{itemize}
\item For every compact set $K\subset\mathbb{C}$, such that $K\cap\overline{\Omega}=\textrm{Ø}$
and $K^{c}$ connected, and every function $h\left(z\right)\in A\left(K\right)$,
there exists a sequence $\left(p_{n},q_{n}\right)_{n=1,2,\ldots}\in\Im$
so that

\textbf{(i).} $f\in\mathfrak{D}_{p_{n},q_{n}}$ for all $n=1,2,\ldots$

\textbf{(ii).} $lim_{n\rightarrow\infty}\left[p_{n}/q_{n}\right]$$_{f}\left(z\right)=h\left(z\right)$
uniformly on $K$ and

\textbf{(iii).} $lim_{n\rightarrow\infty}\left[p_{n}/q_{n}\right]_{f}\left(z\right)=f\left(z\right)$,
uniformly on each compact subset of $\Omega$ .

\end{itemize}
The set of all such $f$ \textquoteright{} s is a dense $G_{\delta}$
subset of $\mathcal{O}\left(\Omega\right)$endowed with the topology
of uniform convergence on compacts. $\blacksquare$

\textbf{Remark 5.6 }In the above theorem it is equivalent to require
convergence of all order derivatives because $\Omega$ and $\left(\overline{\Omega}\right)^{c}$
are open sets (\cite{r17}).

To give a third application of Theorem 5.1, we consider $\Omega\subset\mathbb{C}$
to be a simply connected domain containing $0$, such that $\left\{ \infty\right\} \bigcup\left[\mathbb{C}\setminus\overline{\Omega}\right]$
is connected. We set $L_{k}=\overline{\left\{ z\in\overline{\Omega}:\left|z\right|<n\right\} }$.
The sequence $\left(K_{m}\right)_{m=1,2,\ldots}$ is given by Lemma
5.4. Now the universal functions are smooth on $\partial\Omega$,
in fact they belong to the closure of the set of polynomials in $A^{\infty}(\Omega)$.
We remind that $A^{\infty}(\Omega)$ is the class of all holomorphic
functions $f(z)\in\mathcal{O}(\Omega)$, such that every derivative
$((d^{\left(l\right)}f(z))\lyxmathsym{\textfractionsolidus}(dz^{l}))\;(l=0,1,2,\lyxmathsym{\ldots})$ extends
continuously on $\overline{\Omega}$ (where the closure is taken in
$\mathbb{C}$). The natural topology on $A^{\infty}(\Omega)$ is that
of uniform convergence of all orders derivatives on each compact subset
of $\overline{\Omega}$.

\textbf{Theorem 5.7 }Let $\Omega$ be a simply connected domain of
$\mathbb{C}$ containing $0$, such that $\left\{ \infty\right\} \bigcup\left[\mathbb{C}\setminus\overline{\Omega}\right]$
is connected. Let $\overline{\mathbb{P}}\left(A^{\text{\ensuremath{\infty}}}(\Omega)\right)$denote
the closure of polynomials in $A^{\text{\ensuremath{\infty}}}(\Omega)$.
Let also $\Im\subset\left\{ 0,1,2,\ldots\right\} \times\left\{ 0,1,2,\ldots\right\} $satisfying
condition $\left(\mathcal{H}\right)$. Then there exists a holomorphic
function $f\left(z\right)\in\overline{\mathbb{P}}\left(A^{\text{\ensuremath{\infty}}}(\Omega)\right)$
such that
\begin{itemize}
\item For every compact set $K\subset\mathbb{C}$, such that $K\cap\overline{\Omega}=\emptyset$
and $K^{c}$ connected, and every polynomial $\varphi\left(z\right)$,
there exists a sequence $\left(p_{n},q_{n}\right)_{n=1,2,\ldots}\in\Im$
so that

\textbf{(i).} $f\in\mathfrak{D}_{p_{n},q_{n}}$ for all $n=1,2,\ldots$

\textbf{(ii).} $lim_{n\rightarrow\infty}\frac{d^{l}}{dz^{l}}\left[p_{n}/q_{n}\right]$$_{f}\left(z\right)=\frac{d^{l}}{dz^{l}}\varphi\left(z\right)$,
uniformly on $K$ $\;$($l=0,1,2,...)$and

\textbf{(iii).} $lim_{n\rightarrow\infty}\frac{d^{l}}{dz^{l}}\left[p_{n}/q_{n}\right]$$_{f}\left(z\right)=\frac{d^{l}}{dz^{l}}f\left(z\right)$,
uniformly on each compact subset of $\overline{\Omega}$ .

\end{itemize}
The set of all these functions is a dense $G_{\delta}$ subset of
$\overline{\mathbb{P}}\left(A^{\text{\ensuremath{\infty}}}(\Omega)\right)$.
\textifsymbol[ifgeo]{32}

\textbf{Remark 5.8 }If $\Omega$ is a Jordan domain with rectifiable
boundary, then $\overline{\mathbb{P}}\left(A^{\text{\ensuremath{\infty}}}(\Omega)\right)=A^{\infty}(\Omega)$
(\cite{r38}). \textifsymbol[ifgeo]{80}

Example 3.4 in \cite{r33} section 3 may be transferred to our case.
In this example $\Omega$ is the unit disk. Several properties on
universal Taylor series (: $q=0$) have been established in the literature
especially in the case of the unit disk. One wonder is if they remain
valid in our case. For instance is every universal function non extendable?
Even simpler it is to ask if the radius of convergence of the Taylor
development of a universal function $f$ is exactly equal to the distance
of $0$ from the boundary $\partial\Omega$ of $\Omega$.

In the case $q=0$, it has been examined if we can have universal
Taylor series with respect to several centres simultaneously (\cite{r6}, \cite{r36}, \cite{r37} and \cite{r40}). In our case instead of developing with centre $0$
we can develop with respect to another centre $z_{0}$ and obtain
a formal power series $f\equiv f_{z_{0}}={\displaystyle \sum_{v=0}^{\infty}a_{v}\left(z-z_{0}\right)^{v}}$
with centre $z_{0}$. If $p\in\mathbb{N}_{0}$and $q\in\mathbb{N}_{0}$
are given ($\mathbb{N}_{0}=\left\{ 0,1,2,\text{\ldots}\right\} $),
then \[
\left[p/q\right]_{f,z_{0}}\left(z\right)\]
will denote a rational function\[
\frac{{\textstyle \sum_{v=0}^{p}n_{v}\left(z-z_{0}\right)^{v}}}{{\textstyle \sum_{v=0}^{q}d_{v}\left(z-z_{0}\right)^{v}}},\quad d_{0}=1,\quad n_{p}d_{q}\neq0\]
such that its Taylor development ${\textstyle \sum_{v=0}^{\infty}b_{v}\left(z-z_{0}\right)^{v}}$
with centre $z_{0}$ will have the same $p+q+1$ first coefficients
with $f\equiv f_{z_{0}}$ that is $b_{v}=a_{v}$ for all $v=0,1,\cdots,p+q$.
This rational function may exist or not and it is not necessarily
unique. When such a rational fraction exists is called a Pad\'{e} form
of type $(p,q)$ to the series $f$. A necessary and sufficient condition
for existence and uniqueness is that
\[
\det\left(H^{\left(f\right)}_{q}\left(a_{p-q+1}\right)\right):=\det\left(\underbrace{\begin{array}{ccccc}
a_{p-q+1} & a_{p-q+2} & a_{p-q+3} & \cdots & a_{p}\\
a_{p-q+2} & a_{p-q+3} & a_{p-q+4} & \cdots & a_{p+1}\\
a_{p-q+3} & a_{p-q+4} & a_{p-q+5} & \cdots & a_{p+2}\\
\vdots & \vdots & \vdots &  & \vdots\\
a_{p} & a_{p+1} & a_{p+2} & \cdots & a_{p+q-1}\end{array}}_{q}\right)\neq0.\]

Then we write\[
f\mathfrak{\in D}_{p,q}^{\left(z_{0}\right)}\]
and the Pad\'{e} form is said to be a \textbf{Pad\'{e} approximant of type}
$(p,q)$ to the series $f$.

We are now looking for holomorphic functions $f(z)\in\mathcal{O}(\Omega)$
in a simply connected domain $\Omega$ of $\mathbb{C}$ so that for
every compact set $K\subset\mathbb{C}$, with $K\cap\Omega=\emptyset$
and $K^{c}$ connected, and every function $h(z)\in A(K)$, there
exists a sequence$(p_{n},q_{n})_{n=1,2,\text{\ldots}}\in\Im$ such
that
\begin{itemize}
\item $[p_{n}\lyxmathsym{\textfractionsolidus}q_{n}]_{f,z_{0}}$exists as
a Pad\'{e} form for all $z_{0}\in\Omega$ and all $n=1,2,\lyxmathsym{\ldots}$,
\item $lim_{n\text{\textrightarrow\ensuremath{\infty}}}[p_{n}\lyxmathsym{\textfractionsolidus}q_{n}]_{f,z_{0}}(z)=h(z)$,
uniformly on compact subsets of $(z_{0},z)\in\Omega\times K$ and
\item $lim_{n\text{\textrightarrow\ensuremath{\infty}}}[p_{n}\text{\textfractionsolidus}q_{n}]_{f,z_{0}}(z)=f(z)$,
uniformly on each compact subset of $(z_{0},z)\in\Omega\text{\texttimes}\Omega$.
\end{itemize}

A first question is if there exist universal functions with respect
to all centres. A second question that arises naturally is the following:
does the class of universal functions with respect to one centre coincide
with that with respect to all centres? (See \cite{r37} and \cite{r39}).

We also mention that universal Taylor series do not exist in some
unbounded non-simply connected domains (\cite{r25}). What about universality
of Pad\'{e} approximants?

Finally, we mention that in Theorem 5.1 one can replace the assumption
\textquoteleft{}\textquoteleft{}$K_{m}^{c}$ connected\textquoteright{}\textquoteright{}
by the assumption that \textquoteleft{}\textquoteleft{}$0$ belongs
to the unbounded component of $K_{m}^{c}$\textquoteright{}\textquoteright{},
and then $Q$ will not be anymore a polynomial, but any holomorphic
function in the neighbourhood of $K_{m}$ ($m\in\left\{ 1,2,..\right\} $
being fixed). Then the same result holds provided that $\Im\subset\left\{ 0,1,2,\ldots\right\} \times\left\{ 0,1,2,\ldots\right\} $
satisfies condition $\left(\mathcal{\widehat{H}}\right)$ instead
of $\left(\mathcal{H}\right)$. A difference in the proof is that
$P_{n}$ given by Runge\textquoteright{}s Theorem will not be any
more a polynomial, but it will be a rational function with poles outside
$\Omega\bigcup S$. Then instead of Lemma 2.3 or 2.4, we can use a
variant of Lemma 2.5. In this variant we do not care about the first
coefficients of the Taylor development of the rational function $R$;
thus, the condition $lim_{n\text{\textrightarrow\ensuremath{\infty}}}p_{n}=lim_{n\text{\textrightarrow\ensuremath{\infty}}}q_{n}=+\infty$
is sufficient. Combining the above with the following Lemma, we obtain
Theorem 5.10 bellow.

\textbf{Lemma 5.9 }Let $\Omega$ be a simply connected domain in $\mathbb{C}$
containing $0$. There is a sequence $(\tilde{K}_{m})_{m=1,2,\text{\ldots}}$
of compact subsets of $\Omega^{c}$ where the complement $\tilde{K}_{m}^{c}$
has a finite number of components and $0$ belongs to the unbounded
component of $\tilde{K}_{m}^{c}$, such that the following holds.
\begin{itemize}
\item For every compact set $K\subset\Omega^{c}$ such that $K^{c}$ has
a finite number of components and $0$ belongs to the unbounded component
of $K^{c}$, there exists a $m=1,2,\lyxmathsym{\ldots}$ so that

\textbf{(i).} $K\subset\tilde{K}_{m}$$\;$and

\textbf{(ii).} every component of $K^{c}$ contains a component of
$\tilde{K}_{m}^{c}$.

\end{itemize}
\textbf{Proof }To construct the sequence $(\tilde{K}_{m})_{m=1,2,\text{\ldots}}$,
we start with the sequence $(K_{m})_{m=1,2,\text{\ldots}}$ given
by Lemma 5.2. From each $K_{m}$, we take out a finite union of disjoint
open discs with centers in $\mathbb{Q}+i\mathbb{Q}$ and rational
radii which are included in the interior $int(K_{m})$ of $K_{m}$.
If $int(K_{m})=\mathrm{\slashed{O}}$ then we consider only $K_{m}$
without taking out any disk. Any way even if $int(K_{m})\neq\emptyset$
we also keep $K_{m}$ itself considering that we took off the empty
union of disks. The resulting compact sets are denumerable. An enumeration
gives the family $(\tilde{K}_{m})_{m=1,2,\text{\ldots}}$. One can
easily verify that it has the required property. $\blacksquare$

\textbf{Theorem 5.10 }Let $\Omega$ be a simply connected domain of
$\mathbb{C}$ containing $0$. Let also $\Im\subset\left\{ 0,1,2,\ldots\right\} \times\left\{ 0,1,2,\ldots\right\} $containing
a sequence $(\widetilde{p}_{n},\widetilde{q}_{n})_{n=1,2,\text{\ldots}}$
satisfying $lim_{n\text{\textrightarrow\ensuremath{\infty}}}\widetilde{p}_{n}=lim_{n\text{\textrightarrow\ensuremath{\infty}}}\widetilde{q}_{n}=+\text{\ensuremath{\infty}}$
. Then there exists a holomorphic function $f(z)\in\mathcal{O}(\Omega)$
such that
\begin{itemize}
\item For every compact set $K\subset\mathbb{C}$, such that $K\cap\Omega=\emptyset$
and $K^{c}$ having a finite number of components with $0$ being
in the unbounded component of $K^{c}$, and every function $h(z)$
holomorphic in a neighborhood of $K$, there exists a sequence $(p_{n},q_{n})_{n=1,2,\text{\ldots}}\in\Im$
such that the following hold.

\textbf{(i).} $f\in\mathfrak{D}_{p_{n},q_{n}}$$\bigcap\mathcal{O}(\Omega)$
for all $n=1,2,\ldots$

\textbf{(ii).} $lim_{n\rightarrow\infty}\frac{d^{l}}{dz^{l}}\left[p_{n}/q_{n}\right]$$_{f}\left(z\right)=\frac{d^{l}}{dz^{l}}h\left(z\right)$,
uniformly on $K$ $\;$($l=0,1,2,...)$and

\textbf{(iii).} $lim_{n\rightarrow\infty}\frac{d^{l}}{dz^{l}}\left[p_{n}/q_{n}\right]$$_{f}\left(z\right)=\frac{d^{l}}{dz^{l}}f\left(z\right)$,
uniformly on each compact subset of $\Omega$ $\;$($l=0,1,2,...)$.

\end{itemize}
The set of all these functions is dense and $G_{\delta}$ in $\mathcal{O}(\Omega)$
endowed with the topology of uniform convergence on compacts . $\blacksquare$

\section{THE CASE OF PLANAR DOMAINS WITH ARBITRARY CONNECTIVITY }

Let $\Im\subset\mathbb{N}_{0}^{2}\equiv\left\{ 0,1,2,\text{\ldots}\right\} \times\left\{ 0,1,2,\text{\ldots}\right\} $
containing a sequence $(p_{n},q_{n})_{n=1,2,\text{\ldots}}\in\Im$
with $lim_{n\text{\textrightarrow\ensuremath{\infty}}}p_{n}=lim_{n\text{\textrightarrow\ensuremath{\infty}}}q_{n}=+\text{\ensuremath{\infty}}$
.

Let $\Omega$ be a bounded domain in $\mathbb{C}$ containing $0$.
As usually, $A(\Omega)$ denotes the space of all functions $\phi:\overline{\Omega}\rightarrow\mathbb{C}$
continuous on $\overline{\Omega}$ and holomorphic on $\Omega$. We
suppose that $A(\Omega)$ is endowed with the supremum norm. Further,
$X$ denotes the closure in $A(\Omega)$ of the set of functions holomorphic
in some (varying) neighborhood of $\overline{\Omega}$.

\textbf{Lemma 6.1} Let $\Omega$, $X$ and $\Im$ be as above.Let
also $K\subset\mathbb{C}$ be a compact set with $K\cap\overline{\Omega}=\emptyset.$
Then there exists a $\phi\left(z\right)={\textstyle \sum_{v=0}^{\infty}\varepsilon_{v}z^{v}\in X}$,
such that, for every function $Q(z)$ holomorphic in a neighborhood
of $K$, there exists a sequence $(\widetilde{p}_{n},\widetilde{q}_{n})_{n=1,2,\text{\ldots}}$
in $\Im$, such that

\textbf{(i).} $\phi(z)\in\mathfrak{D}_{\widetilde{p}_{n},\widetilde{q}_{n}}$
for all $n=1,2,\lyxmathsym{\ldots}$ and

\textbf{(ii).} $lim_{n\text{\textrightarrow\ensuremath{\infty}}}[\widetilde{p}_{n}/\widetilde{q}_{n}]_{\phi}(z)=Q(z)$
uniformly on $K$.

The set of all such functions $\phi\in X$  is dense and $G_{\delta}$
in $X$.

\textbf{Proof }For the proof, we may use Baire\textquoteright{}s Theorem.
The essential step is to prove the density of the set
\[
	\bigcup_{\left(p,q\right)\in\Im}
		\left\{
			\begin{array}{c}
				g\left(z\right)\in X:g\in\mathfrak{D}_{p,q}\; and \\
				\sup_{z\in K}\left|\left[p/q\right]_{g}\left(z\right)-Q_{j}\left(z\right)\right|<\frac{1}{s}\;
			\end{array}
		\right\},
\]

\noindent where $Q_{j}$, $j=1,2,\lyxmathsym{\ldots}$ is an enumeration of
rational functions with poles in $(\mathbb{Q}+i\mathbb{Q})\bigcap K^{c}$
such that all coefficients of their numerator and denominator belong
to $\left(\mathbb{Q}+i\mathbb{Q}\right)$. Let $w\left(z\right)$
be a holomorphic function on a neighborhood of $\overline{\Omega}$
and $\varepsilon>0$ be given. Let also $L=\overline{\Omega}$. Using
Runge\textquoteright{}s Theorem, we find two polynomials $A(z)$ and
$B(z)$, with $B(z)\neq0$ on $K\bigcup L$ (in particular $B(0)\neq0$)
and such that \[
\sup_{z\in K}\left|\frac{A\left(z\right)}{B\left(z\right)}-Q_{j}\left(z\right)\right|<\frac{1}{2s}\quad and\quad\sup_{z\in L}\left|\frac{A\left(z\right)}{B\left(z\right)}-w\left(z\right)\right|<\frac{\varepsilon}{2}.\]

Choose a $(p,q)\in\Im$ so that $degA<p$ and $degB<q$. We consider
\[
\widetilde{A}(z)=A(z)+dz^{p}\; and\;\widetilde{B}(z)=B(z)+(cz)^{q},\]
where the constants $d,c\in\mathbb{C}\setminus\left\{ 0\right\} $
will be chosen later on. According to Lemma 2.1, $c$ and $d$ may
be chosen close to zero so that \[
\sup_{z\in K\bigcup L}\left|\frac{\tilde{A}\left(z\right)}{\tilde{B}\left(z\right)}-\frac{A\left(z\right)}{B\left(z\right)}\right|<a,\]

where $0<a<min\left\{ \left(\frac{\varepsilon}{2}\right),\left(\frac{1}{2s}\right)\right\} $.
For $c$ fixed, the Hankel matrix $H{}_{q}^{(A\text{\~\textfractionsolidus}B)}(\delta_{p-q+1})$
of order $q$ at $\delta_{p-q+1}$ for the series\[
\frac{\widetilde{A}\left(z\right)}{\widetilde{B}\left(z\right)}=A\left(z\right)\left[\frac{1}{B\left(0\right)}+\cdots\right]+dz^{p}\left[\frac{1}{B\left(0\right)}+\cdots\right]={\displaystyle \sum_{v=0}^{\infty}\delta_{v}z^{v}}\]

depends linearly on $d$. For $v<p$, the coefficients $\delta_{\nu}$
are independent of $d$. For $v=p$, the coefficient $\delta_{\nu}$
has the form $\delta_{\nu}=d(1\lyxmathsym{\textfractionsolidus}B(0))+c_{p}$,
where $c_{p}$ is independent of $d$ and $(1\lyxmathsym{\textfractionsolidus}B(0))\neq0$.
Thus, following to Lemma 2.2, the determinant $\det\left(H_{q}^{\left(\widetilde{A}/\widetilde{B}\right)}\left(\delta_{p-q+1}\right)\right)$
is a polynomial in $d$ of degree $q$ with leading coefficient $(1\lyxmathsym{\textfractionsolidus}B(0))^{q}(\neq0)$.
The zeros of such a polynomial are finite and we can avoid them by
choosing $d$ close to zero. We infer\[
\frac{\tilde{A}\left(z\right)}{\tilde{B}\left(z\right)}\in\mathfrak{D}_{p,q}\quad and\quad\left[p/q\right]_{\tilde{A}/\tilde{B}}\left(z\right)=\frac{\tilde{A}\left(z\right)}{\tilde{B}\left(z\right)}.\]

Thus, it suffices to set\[
g(z)=\tilde{A}\left(z\right)/\tilde{B}\left(z\right).\]

This gives the result. $\blacksquare$

We now consider an exhausting sequence of compact subsets of $\mathbb{C}\setminus\overline{\Omega}$.
From each member of this sequence, we take out all possible finite
unions of open disks centered at points of $\mathbb{Q}+i\mathbb{Q}$
and with rational radii. Thus, we obtain a sequence $(K_{m})_{m=1,2,\text{\ldots}}$
of compact sets. For any $K$ and every rational function with poles
off $K$, there exists an $m$ so that $K\subset K_{m}$ and the function
is holomorphic in $K_{m}$. Then using Baire\textquoteright{}s Theorem
once more, we obtain the following.

\textbf{Theorem 6.2 }Let $\Im\subset\mathbb{N}_{0}^{2}\equiv\left\{ 0,1,2,...\right\} \times\left\{ 0,1,2,...\right\} $
containing a sequence $(p_{n},q_{n})_{n=1,2,\text{\ldots}}\in\Im$
such that $lim_{n\text{\textrightarrow\ensuremath{\infty}}}p_{n}=lim_{n\text{\textrightarrow\ensuremath{\infty}}}q_{n}=+\text{\ensuremath{\infty}}$
. Let $\Omega$ be a bounded domain of $\mathbb{C}$ containing $0$.
Let also $X$ denotes the closure in $A(\Omega)$ of the set of functions
holomorphic in a (varying) neighborhood of $\overline{\Omega}$. Then
there exists a $\phi(z)\in X$, such that the following holds.
\begin{itemize}
\item For all compact sets $K\subset\mathbb{C}$, such that $K\cap\overline{\Omega}=\emptyset$
and all functions $Q(z)$ holomorphic in a neighborhood of $K$ there
exists a sequence $(p_{n},q_{n})_{n=1,2,\text{\ldots}}\in\Im$ so
that

\textbf{(i).} $\phi\in\mathfrak{D}_{p_{n},q_{n}}$ for all $n=1,2,\lyxmathsym{\ldots}$
and

\textbf{(ii)}. $lim_{n\text{\textrightarrow\ensuremath{\infty}}}[p_{n}/q_{n}]_{\phi}(z)=Q(z)$
uniformly on $K$.

\end{itemize}
The set of all such $\phi$\textquoteright{}s is a dense$G_{\delta}$
subset of $X$. $\lyxmathsym{\textifsymbol[ifgeo]{80}}$

\textbf{Remark 6.3 }In Theorem 6.2 one can also require

\textbf{(iii).} $lim_{n\text{\textrightarrow+\ensuremath{\infty}}}[p_{n}/q_{n}]_{\phi}(z)=\phi(z)$
in $X.$

The reason is that in the proof of Lemma 6.1 we automatically have
$\phi-[p\lyxmathsym{\textfractionsolidus}q]_{\phi}\equiv0$ because
$\phi\left(z\right)=\widetilde{A}\left(z\right)/\widetilde{B}\left(z\right)\in\mathfrak{D}_{p_{n},q_{n}}$
is a rational function. We see therefore that it is a generic property
of holomorphic functions in a domain to be limits of some Pad\'{e} approximants
of them; see also (\cite{r8}).

\textbf{Remark 6.4 }With some more effort one can prove a version
of Theorem 6.2 where $X$ is replaced by $X^{\infty}$ the closure
in $A^{\infty}(\Omega)$ of the set of functions holomorphic in some
(varying) neighborhood of $\overline{\Omega}$. The set $A^{\infty}(\Omega$)
is the space of all holomorphic functions $f(z)$ in $\Omega$ such
that every derivative $f^{(l)}(z)$ ($l=0,1,2,\lyxmathsym{\ldots}$)
extends continuously on $\overline{\Omega}$. The topology of $A^{\infty}(\Omega)$
is defined by the sequence of seminorms $\left\Vert f\right\Vert _{l}=\sup_{z\in\overline{\Omega}}\left|f^{\left(l\right)}\left(z\right)\right|$,
provided $\Omega$ is bounded.

\textbf{Remark 6.5 }The assumption in Theorem 6.2 that $\Omega$ is
bounded is not essential. In fact, if $L=\overline{\left\{ z\in\Omega:\left|z\right|<n\right\} }$
is a compact subset of $\overline{\Omega}$ and a pole $z_{0}$ of
the approximating rational function belongs to $\Omega\setminus L$,
then $z_{0}$ lies in the same component of $(\mathbb{C}\bigcup\left\{ \infty\right\} )\setminus(L\bigcup K)$
with some point $z_{1}\in(\mathbb{C}\bigcup\left\{ \infty\right\} )\setminus(\overline{\Omega}\bigcup K)$.
This is possible since $dist(\overline{\Omega},K)>0$. Thus, the approximating
rational function is finite on $\overline{\Omega}$.

\end{document}